\documentclass[12pt]{article}

\title{Diophantine Approximation on varieties IV: Derivated algebraic distance and derivative metric B\'ezout Theorem}

\author{Heinrich Massold}

\usepackage{bezier,amsfonts}

\usepackage{amsmath,amsfonts,bbm,amssymb}

\newtheorem{Satz}{}[section]

\newcommand{\satz}[1]{\vspace{2mm} \begin{Satz}{\bf #1}}

\newcommand{\proof}{\vspace{4mm} {\sc Proof}\hspace{0.2cm}}

\newcommand{\la}{\langle}

\newcommand{\ra}{\rangle}

\newcommand{\di}{{\mbox{div}}}

\newcommand{\sgn}{\mbox{sgn}}

\newcommand{\R}{{\mathbbm R}}

\newcommand{\Pe}{{\mathbbm P}}

\newcommand{\Z}{\mathbbm{Z}}

\newcommand{\N}{\mathbbm{N}}

\newcommand{\CC}{{\cal C}}

\newcommand{\CX}{{\cal X}}

\newcommand{\CY}{{\cal Y}}

\newcommand{\CZ}{{\cal Z}}

\newcommand{\CL}{{\cal L}}

\newcommand{\C}{{\mathbbm C}}

\newcommand{\A}{{\mathbbm A}}

\newcommand{\spec}{\mbox{Spec}}

\begin{document}

\parindent0mm

\maketitle

\thispagestyle{empty}

\tableofcontents

\section{Introduction}

Let $\Pe^t = \Pe(\Z^{t+1})$ be the $t$-dimensional projective space
over $\spec \Z$, $\Pe(W) \subset \Pe^t_\C$ a projective subspace of 
codimension $q$, and $X \in Z_{eff}(\Pe^t_\C)$ an effective cycle.
We say that $X$ is reguar with respect to $\Pe(W)$ if the 
irreducible components of codimension $< t-q$ intersect $\Pe(W)$
properly, and the irreducible components of codimension $\geq t-q$ do
not meet $\Pe(W)$. For regular cycles the algebraic distance
\[ D(Y,\Pe(W)) \in \R \]
is defined in \cite{App1}, section 4.1.
Further for $p+q \leq t+1$, and effective cycles $X,Y$
of pure codimensions $p,q$ that intersect properly, the algebraic
distance
\[ D(X,Y) \in \R \]
is defined in \cite{App1}, Definition 4.1.
For $x,y$ points in $\Pe^t(\C)$ denote by $|x,y|$ their Fubini-Study distance,
i.\@ e.\@ $\sin (x,y)$. The logarithm of the distance is a nonpositive number.

There are the following Theorems for the algebraic distance. 

\satz{Theorem} \label{arbezout}
For properly intersecting cyles $\CX, \CY \in Z_{eff}(\Pe^N_\Z)$ with
base extensions $X,Y$ to $\C$ the equality
\[ D(Y,Z) = h(\CY . \CZ) - \deg Z h(\CY) - \deg Y h(\CZ) 
   - \sigma_t \log 2 \deg Y \deg Z \]
holds.
\end{Satz}

\proof
\cite{App1}, Scholie 4.3.

\satz{Theorem} \label{bezout}
With the previous Definition, let $\CX,\CY$ be effective cycles 
intersecting properly, and $\theta$ a point in 
$\Pe^N(\C) \setminus (supp(X_\C\cup Y_\C))$. 
\begin{enumerate}

\item
There are effectively computable constants $c,c'$ only depending
on $t$ and the codimenion of $X$ such that

\[ \deg(X) \log |\theta,X(\C)| \leq D(\theta,X) + c \deg X \leq 
   \log |\theta,X(\C)| + c' \deg X, \]
\item
If $\CX = \di(f)$ is an effective cycle of codimension one,
\[ h(\CX) \leq \log |f_D|_{L^2} + D \sigma_t, \quad \mbox{and} \]
\[ D(\theta,X) + h(\CX) = \log |\la f| \theta \ra| + D \sigma_{N-1}, \]
where the $\sigma_i's$ are certain constants, and
$|\la f|  \theta \ra|$ is taken to be the norm of the evaluation of 
$f \in Sym^D(E) = \Gamma(\Pe^N,O(D))$ 
at a vector of length one representing $\theta$.

\item 
For $p+q \leq N+1$,
assume that $\CX$, and $\CY$ have pure codimension $p$, and $q$ 
respectively.
There exists an effectively computable positive constant $d$,
only depending on $t$, and a map 
\[ f_{X,Y}: I \to \underline{\deg X} \times \underline{\deg Y} \]
from the unit interval $I$ to the set of natural numbers
less or equal $\deg X$ times the set of natural numers less or equal
$\deg Y$ such that $f_{X,Y}(0) = (0,0), f_{X,Y}(1) = (\deg X, \deg Y)$,
and the maps $pr_1 \circ f_{X,Y}:I \to \underline{\deg X}, 
pr_2\circ f_{X,Y}: I \to \underline{\deg Y}$
are monotonously increasing, and surjective, fullfilling: For every
$T \in I$, and $(\nu,\kappa) = f_{X,Y}(T)$, the inequality
\[ \nu \kappa \log |\theta,X+Y| + D(\theta,X.Y) + h(\CX . \CY) \leq \]
\[  \kappa D(\theta,X) + \nu D(\theta,Y) +
   \deg Y h(\CX) + \deg X h(\CY) + d \deg X \deg Y \]
holds. 
\item
In the situation of 3, if further
$|\theta,X+ Y| = |\theta,X|$, then
\[ D(\theta,X . Y) + h(\CX . \CY) \leq D(\theta,Y) +
   \deg Y h(\CX) + \deg X h(\CY) + d' \deg X \deg Y \]
with $d'$ a constant only depending on $N$.
\end{enumerate}
\end{Satz}

\proof
\cite{App1} Theorem 2.2.

\vspace{2mm}

Let now $\Pe(W) \subset \Pe^t(\C)$ be a subspace of dimension $q$, and
$\partial^I$ a real differential operator on the Grassmannian $G_{t+1-q,t+1}$
of $t+1-q$-dimensional subspaces of $\C^{t+1}$ with respect to some affine
chart, where $I=(i_1, \ldots, i_{2q(t+1-q)}$ is a multiindex of order
$|I|= i_1 + \cdots + i_{2q(t+1-q)}$. (More details will be given later on).
Defines the derivated algebraic distance of order $S$ of $X$ to $\Pe(W)$ as
\[ D^S(Y,\Pe(W))) := 
   \sup_{|I| \leq S}  
   \log |\partial^I (\exp D(X,\Pe(W)))|. \]
There are the following Theorems for this derivated algebraic distance.

\satz{Theorem} \label{hypeb}
For $s,D \in \N$, and $f \in \Gamma(\Pe^t,O(D))$ let $F$ be the polynomial
of degree at most $D$ in $t$ variables that corresponds to $f$ with respect
to affine coordinates of $\Pe^t$ centered at $\theta$. Then,
with some positive constant $c$ only depending on $t$, 
\[ D^S(\di f,\theta) \leq \sup_{s \leq S |J| = s} 
  \log \left| \left( \frac{\partial^s}{(\partial z_1)^{j_1} \cdots
   (\partial z_t)^{j_t}} f \right) (0) \right| + c (s+D) \log (SD). \]
\end{Satz}

in the following Theorems the $O$-notation always signifies that the
respective inequalities hold modulo a fixed contant only depending on $t$ and
codimnstions of cylces times the term inside the $O$-bracket.

\satz{Theorem} \label{main2}
Let $Z$ be an effective cycle of pure codimension $p$ in $\Pe^t$,
and $\theta$ a point not contained in the support of $Z$, and let
$|\cdot, \cdot|$ denote the Fubini-Study distance in $\Pe^t$.

There is a projective subspace $\Pe(F) \subset \Pe^t_\C$ of codimension
$t-p$ intersecting $Z$ properly such that with
$z_1, \ldots, z_{\deg Z}$ the points in the intersection $\Pe(F).Z$
counted with multiplicity,
$|z_1,\theta| \leq \cdots \leq |z_{\deg Z},\theta|$, and for $S < \deg Z/3$
the equalities
\[ D^S(Z,\theta) = \sum_{i=S+1}^{\deg Z} \log |z_i,\theta| + 
   O(S \log \deg Z), \]
\[ 2\sum_{i=S+1}^{\deg Z} \log |z_i,\theta| \leq
   D^{3S}(Z,\theta) + O((\deg Z+ S) \log (S \deg Z)) \]
hold.
\end{Satz}

\satz{Corollary}
The derivated algebraic distance is a negative number modulo \\
$O((\deg Z+ S) \log (S \deg Z))$.
\end{Satz}

Next, for $n \in \N$ denote $\underline{n}$ the set of natural numbers less
or equal $n$ including $0$, and for $Z_0,Z_1$ effective cyles, let 
$f =(f_0,f_1): \; \underline{\deg Z_0+\deg Z_1} \to 
\underline{\deg Z_0} \times \underline{\deg Z_1}$ be a a path from 
$(1,1)$ to $(\deg Z_0, \deg Z_1)$ such that in each step exactly
one of the coordinates increases. If in the $k$th step the coordinate
$i$ increases, set $i_k = i$. 

\satz{Theorem} \label{DMBT1}
For any effective cycles $Z_0, Z_1 \in Z_{eff}(\Pe^t)$ that intersect properly,
and $\theta \in \Pe^t_\C$ a point not contained in the support
of $Z_0 . Z_1$, there is a path $f$ such that
\[ 2D(Z_0,Z_1) + 2D(Z_0.z_1, \theta) \leq \sum_{k=1}^n D^{f_{i_k}(k)}
   (Z_{i_k},\theta)+ \]
\[ O((\deg Z_0 \deg Z_1+ S) \log (S \deg Z_0 \deg Z_1)). \]
\end{Satz}

\satz{Corollary} \label{cor1}
With any $l \leq \deg Z_0 + \deg Z_1$, and $(\nu_0,\nu_1) = f(l)$,
\[ 2D(Z_0,Z_1) + 2D(Z_0.Z_1,\theta) \leq \]
\[ \nu_0  D^{3\nu_1}(Z_1,\theta) + 
   \nu_1  D^{3\nu_0}(Z_0,\theta) + 
   O((\deg Z_0 \deg Z_1+ S) \log (S \deg Z_0 \deg Z_1)). \]
\end{Satz}

This immediately implies

\satz{Corollary} \label{cor2}
For every $S \leq \deg Z_1/3$, 
\[ 2D(Z_0,Z_1) + 2D(Z_0.Z_1, \theta) \leq 
   \mbox{max}(S D(Z_0,\theta), D^{3S-1}(Z_1,\theta))+ O(\deg X \deg Y). \]
\end{Satz}

\satz{Theorem} \label{DMBT2}
For any $S \leq \deg X \deg Y/9$, there are natural numbers
$\nu_0,\nu_1$ with $\nu_0 \nu_1 \leq S$, and a path $f$ such that
and a function $h_S: \underline{\deg Z_0 + \deg Z_1} \to \N$ with
$h_S(k) = 0$ for $k \geq k_0$ such that 
\begin{enumerate}

\item
\[ 2D(X,Y) + 2D^S(X.Y,\theta) \leq \sum_{k=1}^n
   D^{3(f_{i_k}(k)-h_S(k))} (Z_{i_k},\theta). \]

\item
For any $k \leq \deg Z_0 + \deg Z_1$ greater or equal $k_0$,
and $(\bar{\nu}_0,\bar{\nu}_1) = f(k)$,
\[ 2(\bar{\nu}_0-\nu_0) (\bar{\nu}_1-\nu_1) 
   \log |Z_0+Z_1,\theta| + 2D^S(Z_0.Z_1,\theta)+ D(Z_0,Z_1) \leq \]
\[ (\bar{\nu}_0 - \nu_0) D^{3\kappa_0}(Z_1,\theta) + 
   (\bar{\nu}_1-\nu_1) D^{3\kappa_0}(Z_0,\theta) + \]
\[ O((\deg Z_0 \deg Z_1+ S) \log (S \deg Z_0 \deg Z_1)). \]
\end{enumerate}
\end{Satz}

\satz{Corollary} \label{cor3}
For any $S \leq \deg Z_0 \deg Z_1/9$, there are numbers $\nu_0 \nu_1$ with
$\nu_0 \nu_1 \leq S$ such that for any $l \leq \deg Z_0 + \deg Z_1$ greater
or equal $k_0$ with 
$(\bar{\nu}_0, \bar{\nu}_1) = f(l)$,
\[ 2D(Z_0,Z_1) + 2D^S(Z_0.Z_1,\theta) \leq \]
\[ (\bar{\nu}_1-\nu_1)D^{3\bar{\nu}_0}(Z_0,\theta) + 
   (\bar{\nu}-\nu_0)D^{3\bar{\nu}_1}
   (Z_1,\theta) + O((\deg Z_0 \deg Z_1+ S) \log (S \deg Z_0 \deg Z_1)). \]
\end{Satz}

\satz{Corollary} \label{cor4}
Let $S_0 \leq \deg Z_0/3$, $S_1 \leq \deg Z_1/3$ be natural numbers, and
$S= S_0 S_1$. Then,
\begin{eqnarray*}
2D(Z_0,Z_1) + 2 D^S(Z_0 . Z_1, \theta) &\leq &
   \mbox{max} (S_1 D^{9S_0}(Z_0,\theta), S_0 D^{9 S_1}(Z_1,\theta)) + \\
&& O((\deg Z_0 \deg Z_1+ S) \log (S \deg Z_0 \deg Z_1)). 
\end{eqnarray*}
\end{Satz}

\vspace{2mm}

{\bf Remarks:} 1. I strongly conjecture that Theorem \ref{main2} as well
as Theorems \ref{DMBT1}, \ref{DMBT2}, and their corollaries still hold if
the factor $2$ before $D(\theta,X.Y)$ and $D^S(\theta,X.Y)$ is dropped, and
possibly also if the $3$ in the exponent on the right hand side is replaced by 
some smaller number greater or equal $1$. In order to obtain this, one would
only have to improve Lemma \ref{hilfzwei} in this respect; however, I don't
know right now how to do that. For the applications of the Theorems and
Corollaries to Diophantine Approximation and algebraic independence
theory this improvement would be insubstantial.

\vspace{2mm}

2. Throughout this paper, constants entailed by the notaion $O(\cdots)$
always depend only on $t$, and the dimensions of cycles involved
in the context. As, there are always only finitely many cycles involved,
the constants can also be assumed to be depending only on $t$.

\vspace{2mm}

3. To my knowledge, in the literature, one special case of 
the above Theorems and Corollarys is known, namely Corollary \ref{cor2} in the
case $t=1$, $\mbox{codim}Z_0 = \mbox{codim} Z_1 =1$.
See \cite{LR}, preuve du corollaire 3.

\vspace{2mm}


Recall that in \cite{App1}, section 4 there were given 3 alternative 
definitions of the algebraic distance. The algebraic distance
$D_\infty(\theta, \cdot)$ is not additive on the cycle group, and has
some other deficencies; therefore it is probably not possible to prove
the derivative metric B\'ezout for $D_\infty$. Proofs will be given for
$D_{Ch}$ and $D_1$. 
\vspace{2mm}

This paper heavily depends on part one (\cite{App1}) of this series
on diophantine approximation on varieties.
and can possible not be read independentyl of it. It does not however 
presuppose any knowledge of part 2 and part 3.

\section{Sharp decomposition of the algebraic distance}

Recall the following notations from \cite{App1}. If $G=G_{q,t}$ is the
Grassmannian of $q$-dimensional subspaces of $\C^{t+1}$, then for
a subspace $\Pe(W) \subset \Pe^t$ of dimension $r \leq q$ the sub Grassmannian 
of $G$ consisting of the spaces that contain $W$ is denoted $G_W$, and
for $\Pe(F) \subset \Pe^t$ a subspace of dimension $p \geq q$ the
sub Grassmannian of spaces being contained in $F$ is denoted $G^F$.

Let $\Pe(F) \subset \Pe_\C$ be a subspace of codimension $r$, and $\pi$ the map
\[ \pi : \Pe^t \setminus \Pe(F^\bot) \to \Pe(F), \quad
   [v,w] \mapsto [v], \quad v \in F, w \in F^\bot. \]
For any sub variety $X \subset \Pe(F)$ of codimension $p$, the
closure $X_F := \overline{\pi^{-1}(X)}$ is a subvariety of codimension
$p$ in $\Pe^t$ with the same degree as $X_F$. This induces a map
$\pi^*:Z^p(\Pe(F)) \to Z^p(\Pe^t), X \mapsto X_F$ with left inverse
$X \mapsto X . \Pe(F)$.
For two effective cycles $X \in Z^p(\Pe(F)), Y \in Z^q(\Pe(F))$, denote by
$D_{\bullet}^{\Pe(F)}(X_F,Y_F)$ their algebraic distance  as cycles in
$\Pe(F)$.

\vspace{2mm}

In \cite{App1}, Theorem 4.11, and Proposition 4.16,
for $Z$ an effective cycle in $\Pe^t_\C$, of codimension $p$
and $\Pe(F) \supset \Pe(W)$ subspaces of codimensions $r \geq t-p, q > r$
respectively, regular with respect to $X$ the relations
\[ D(\Pe(W),Z) \leq D(\Pe(W),\Pe(F).Z) + c_1 \deg Z \leq
   D(\Pe(W),Z) - D(\Pe(F),Z) + c_2 \deg Z \]
with $c_1,c_2$ constants depending only on $p,q,r$, and $t$,
were proved, and thereby the algebraic distance of $Z$ to $\Pe(W)$ modulo
$O(\deg Z)$ is reduced to the algebraic distance of $\Pe(W)$ to 
$Z . \Pe(F)$ and the algebraic distance of $\Pe(F)$ ot $Z$. 

If one wants to consider derivatives of algebraic distances, this decomposition
is not good enough, because the derivatives of two functions may have an
arbitrarily big difference even if there values don't differ very much. One
needs the following sharper decomposition.

\satz{Proposition} \label{raumraum}
With the above notations, there is a positive constant $c_1$ only depending
on $p,q,r$, and $t$ such that
\[  D(\Pe(W),Z) + c_1 \deg Z = D(\Pe(W),Z.\Pe(F)) + D(Z,\Pe(F)). \]
\end{Satz}

The proof will use two Lemmas.

\satz{Lemma} \label{raumraum1}
Let $p+r \leq t$, $X \in Z^p_{eff}(\Pe_\C)$, and $\Pe(F) \subset \Pe(\C)$
a subspace of codimension $r$ that intersects $X$ properly. Let further 
$\Pe(W) \subset \Pe(F)$ be a subspace of codimension $q \geq r$ that is
regular with respect to $X$. Then
\[ D(Z_F,\Pe(W)) - D(Z,\Pe(W)) = D(Z_F,\Pe(F)) - D(Z,\Pe(F). \]
\end{Satz}

\proof
If $q \leq t+1-p$, this is \cite{App1}, Proposition 5.1.

If $q > t+1-p$, the Lemma will be proved for $D_{Ch}$, and $D_G$ successively,
firstly for $D_{Ch}$:
Let $\Pe(V) \subset \Pe^t_\C$ be a subspace of codimension $t+1-p$ that does
not intersect $Z$, and fullfills $\Pe(W) \subset \Pe(V) \subset \Pe(F)$.
By \cite{App1}, Proposition 5.1, and Proposition 4.14.2,
\[ D(Z_F,\Pe(F)) - D(Z,\Pe(F)) = D(Z_F, \Pe(V))- D(Z, \Pe(V)) = \]
\begin{equation} \label{raumrauma} 
D_{Ch}(Z_F, \Pe(V)) - D_{Ch}(Z,\Pe(V)) .
\end{equation}

Next, we repeat the construction of the cycle deformation in \cite{App1}, 
section 5. For $\lambda \in \C$, 
\begin{equation} \label{deform}
(F, \pi, \psi_\lambda, X_\lambda) 
\end{equation}
is defined as follows.
Let $X$ be a subvariety of codimension $p$ in $\Pe^t$, further
$F \subset \C^{t+1}$ a sub vector space that is regular with respect to $X$,
and $F^{\bot} \subset \C^{t+1}$
the orthogonal complement of $F_\C$ with respect to the canonical inner 
product on $\C^{t+1}$, and $\Pe(F_\C^{\bot})$ the corresponding projective
subspace of $\Pe^t$. Consider the map from above
\[ \pi: \Pe(E_\C) \setminus \Pe(F_\C^\bot) \to \Pe(F_\C), \quad
   [v \oplus w] \mapsto [v].  \]
For each $\lambda \in \C^*$, there is the automorphism
\[ \psi_\lambda: 
   \Pe(E_\C) \to \Pe(E_\C), \quad [v \oplus w] \mapsto [\lambda v + w]. \]
For any effective cycle $X$ in $\Pe^t$ that intersects $\Pe(F)$ properly,
define, $\Phi$ as the subvariety of $(\Pe^t_\C)^{t+1-p} \times \A_\C$ 
given as the Zariski closure of the set
\begin{equation} \label{abschl}
\{ (\psi_\lambda(x),\ldots,\psi_\lambda(x),\lambda) 
   \in (\Pe^t_\C)^{t+1-p} \times \C^*| x \in X(\C). \}.
\end{equation}
Then, $\Phi$ intersects $\Pe(F) \times \A^1$ properly.
Further, for $\lambda \in \C$, and $y$ the coordinate of the affine line,
the divisor $\Phi_\lambda$ corresponding to
the restriction of the function $y - \lambda$ to $\Phi$ is 
a proper intersection of $\Phi$ and the zero set of $y-\lambda$ and is
of the form $Z_\lambda^{t+1-p} \times \{\lambda\}$, for some subvariety 
$Z_\lambda$ of $\Pe^t_\C$, and for $\lambda \neq 0$, we have
$Z_\lambda = \psi_\lambda(X)$. The specialization $\Phi_0$ equals
\[ \Phi_0 = \overline{\pi^*(X . \Pe(F)} = 
   \overline{\pi^*(X_{\lambda} . \Pe(F))}, \]
for arbitrary $\lambda \in \C^*$. (See  \cite{BGS}, p.\@994.)

Recall the correspondence from \cite{App1}, section 5,

\special{em:linewidth 0.4pt}
\linethickness{0.4pt} 
\vspace{4mm}
\hspace{35mm}
\begin{picture}(200.00,50.00)
\put(65.00,35.00){\vector(-4,-3){28.00}}
\put(100.00,45.00){\makebox(0,0)[cc]{$\bar{C}=C \times \A^1$}}
\put(15.00,5.00){\makebox(0,0)[cc]{$(\Pe^t)^{t+1-p} \times \A^1$}}
\put(180.00,5.00){\makebox(0,0)[cc]{$(\check{\Pe}^t)^{t+1-p} \times \A^1$}}
\put(40,28){\makebox(0,0)[cc]{$F$}}
\put(152,28){\makebox(0,0)[cc]{$G$}}
\put(130.00,35.00){\vector(4,-3){28.00}}
\end{picture}

\vspace{3mm}

and the map $\Delta : \Pe^t \times \A^1 \to (\Pe^t)^{t+1-p} \times \A^1$ 
from the proof of \cite{App1}, Proposition 4.16.3:
$C$ is the Chow correspondence from (\ref{chow}); and $F,G,\Delta$ are defined
by taking the identity on the second factor.

With $(F, \pi, \psi_\lambda, Z_\lambda)$ 
for an effective cycle $Z$ of codimension as above, and the corresponding 
cycle $\Phi \subset \Pe^t \times \A^1$, the cycle 
$\bar{\Phi} := G_* F^* \Delta_* \Phi$ intersects 
$(\Pe(\check{W}))^{t+1-p} \times \A^1$ 
properly, and the intersection of $\bar{\Phi}$ with 
$(\check{\Pe}^t)^{t+1-p} \times \{y\}$ is likewise proper.
For any $\lambda \in \C$ the cycle 
$G_* F^* \Delta_* (Z_\lambda \times \{ \lambda \})$ equals
$Ch(Z_\lambda) \times \{ \lambda \}$. Take $g_{\Pe(\check{W}^p)}$
the normalized Green form of $\Pe(\check{W})^{t+1-p}$ in 
$(\check{\Pe}^t)^{t+1-p}$, and define
\[ \varphi_W(\lambda) = \int_{(\check{\Pe}^t)^{t+1-p}} 
   \delta_{Ch(Z_\lambda)} g_{\Pe(\check{W})^{t+1-p}} 
   \bar{\mu}^{(q-1)(t+1-p)}, \]
\[ \varphi_V(\lambda) = \int_{(\check{\Pe}^t)^{t+1-p}} 
   \delta_{Ch(Z_\lambda)} g_{\Pe(\check{V})^{t+1-p}} 
   \bar{\mu}^{(t-p)(t+1-p)}, \]
where $\bar{\mu} = c_1(\overline{O(1, \ldots,1)})$ (see \cite{App1}, 
section 3.3 for details). By definition,
\begin{equation} \label{raumraumb}
D_{Ch}(\Pe(W),Z_\lambda) = \frac1{{(t+1-p)(t-p) \choose t-p, \ldots, t-p}}
\varphi_W(\lambda), \quad
D_{Ch}(\Pe(V),Z_\lambda) = \frac1{{(t+1-p)(r-1) \choose r-1, \ldots, r-1}}
\varphi_V(\lambda),
\end{equation}
where ${(t+1-p)(r-1) \choose r-1, \ldots, r-1}$ is the multinomial coefficient
$\frac{((t+1-p)(r-1))!}{(r-1)!^{t+1-p}}$.
From the proof of \cite{App1}, Proposition 4.16.3, it is clear that
there are smooth real functions $\chi_W, \chi_F: \R \to \R$ such that
$\varphi_W(\lambda) = \chi_W(\log |\lambda|)$, 
$\varphi_V(\lambda) = \chi_V(\log |\lambda|)$, and further that
$\chi_W'(0) = \chi_V' = 0$, and $\chi_W'', \chi_V''$ are nonnegative; 
consequently $\varphi_W(\lambda) \geq \varphi_W(0)$,
$\varphi_V(\lambda) \geq \varphi_V(0)$. Finally, with
$\mu_i$ being the Fubini-Study form on the $i$th factor of $(\Pe{t})^{t+1-p}$,
\[ d d^c([\varphi_W]) = 
   {pr_2}_*(\delta_{\bar{\Phi}} pr_1^* (\mu_1^{t+1-q} \cdots 
   \mu_{t+1-p}^{t+1-q} \; \bar{\mu}^{q(t+1-p)-1})) = \]
\[ {(t+1-p)(t-p) \choose t-p,\ldots,t-p} pr_{2*} \left( \delta_{\bar{\Phi}}
   pr^*_1(\mu_1^t \cdots \mu_{t+1-p}^t) \right),  \]
and similarly
\[d d^c([\varphi_V]) = 
  {(t+1-p)(r-1) \choose r-1,\ldots,r-1} pr_{2*} \left( \delta_{\bar{\Phi}}
  pr^*_1(\mu_1^t \cdots m_{t+1-p}^t) \right), \]
hence
\begin{equation} \label{raumraumc}
\frac1{{(t+1-p)(t-p) \choose t-p,\ldots,t-p}} d d^c ([\varphi_W]) =
\frac1{{(t+1-p)(r-1) \choose r-1,\ldots,r-1}} d d^c ([\varphi_V]) =
pr_{2*} \left( \delta_{\bar{\Phi}} pr^*_1(\mu_1^t \cdots m_{t+1-p}^t) \right),
\end{equation}
which does not depend on $\Pe(W)$, and $\Pe(V)$ but only on the numbers 
$p$ and $t$, the subspace $\Pe(F)$, and the cycle $Z$. 
We get
\[ \frac1{{(t+1-p)(t-p) \choose t-p,\ldots,t-p}} \chi_W'' = 
   \frac1{{(t+1-p)(r-1) \choose r-1,\ldots,r-1}}\chi_V'', \]
consequently, since $\chi_W'(0) = \chi_V(0)'=0$,
\begin{equation} \label{abl1}
\frac1{{(t+1-p)(t-p) \choose t-p,\ldots,t-p}}
(\varphi_W (0) - \varphi_W(\lambda) = 
\frac1{{(t+1-p)(r-1) \choose r-1,\ldots,r-1}}
(\varphi_F (0) - \varphi_F(\lambda)). 
\end{equation}
Together with (\ref{raumraumb}), this implies
\[ D(\Pe(W),Z_F) - D(\Pe(W,Z) = D(\Pe(W),Z_0)- D(\Pe(W),Z_1) = \]
\[ D(\Pe(V),Z_F)- D(\Pe(V),Z), \]
which together with (\ref{raumrauma}) entails the Lemma.

\proof {\sc of Proposition \ref{raumraum} for $D_G$:}
Again, if $\Pe(V)$ is a subspace of codimension $t+1-p$ interesecting $Z$
properly such that $\Pe(V) \subset \Pe(F)$, then, by the proof of
\cite{App1}, Proposistion 4.14,
\[ D(\Pe(F), Z_F)- D(\Pe(F),Z) = D(\Pe(V),Z_F) - D(\Pe(V),Z) = \]
\[ D_G(\Pe(V,Z_F))- D_G(\Pe(V,Z)). \]
With $(F,\pi,\psi_\lambda, X_\lambda)$ as in (\ref{deform}), and the
correspondence

\special{em:linewidth 0.4pt}
\linethickness{0.4pt} 
\vspace{4mm}
\hspace{35mm}
\begin{picture}(200.00,50.00)
\put(65.00,35.00){\vector(-4,-3){28.00}}
\put(100.00,45.00){\makebox(0,0)[cc]{$\bar{C}=C \times \A^1$}}
\put(15.00,5.00){\makebox(0,0)[cc]{$\Pe^t \times \A^1$}}
\put(180.00,5.00){\makebox(0,0)[cc]{$G \times \A^1$}}
\put(40,28){\makebox(0,0)[cc]{$F$}}
\put(152,28){\makebox(0,0)[cc]{$G$}}
\put(130.00,35.00){\vector(4,-3){28.00}}
\end{picture}
\vspace{-11mm}
\begin{equation} \label{grascor} \end{equation}

\vspace{4mm}

where again $F$, and $G$ are defined by taking the identity on the
second factor, and the intersections of $\bar{\Phi}$ with
$G \times \{y\}$ are proper. Thus, we can proceed just in the
case $D_{Ch}$, and define
\[ \varphi_W(\lambda) = \int_G \delta_{V_{Z_\lambda}}
   g_{G_W} \mu_G^{(t+1-p)(t+1-q)}, \quad 
   \varphi_V(\lambda) = \int_G \delta_{V_{Z_\lambda}} g_{G_V}, \]
where $\mu_G = c_1(L_G)$ with $L_G$ the canonical line bundle on $G$.
(see \cite{App1}, section 3.2.) We have again
\[ [\varphi_W] = (pr_2)_* (\delta_\Phi pr_1^* (g_{G_W} 
   \mu_G^{(t+1-p)(t+1-q)})),
   \quad [\varphi_V] = {pr_2}_* (\delta_\Phi pr_1^* (g_{G_V})), \]
and the calculations
\[ d d^c (\delta_{\bar{\Phi}} pr_1^* (g_{G_W})) +
   \delta_{\bar{\Phi} . pr_1^*(G_W)} = \delta_{\bar{\Phi}}
   pr_1^* (\omega(g_{G_W})), \]
\[ d d^c (\delta_{\bar{\Phi}} pr_1^* (g_{G_V}) +
   \delta_{\bar{\Phi} . pr_1^*(G_V)} = \delta_{\bar{\Phi}}
   pr_1^* (\omega(g_{G_V})) \]
imply this time
\[ d d^c([\varphi_W]) = 
   pr_{2*} (\delta_{\bar{\Phi}} pr^*_1 (\omega(g_W) \mu_G^{(t+1-p)(t+1-q)})),\]
\[ d d^c([\varphi_V]) = 
   pr_{2*} (\delta_{\bar{\Phi}} pr^*_1 (\omega(g_V)). \]
Now, $G_V$ is the single point $V$ in $G$, hence $\omega(g_V)$ is the
canonical generator of the one dimensional space of harmonic forms
$H^{p(t+1-p),p(t+1-p)}(G)$. Further, by the intersection theory on $G$
the form $\omega(g_W) \mu_G^{(t+1-p)(t+1-q)}$ likewise equals this generator. 
We get
\[ d d^c ([\varphi_W]) = d d^c ([\varphi_V]), \]
and the rest of the proof is in complete analogy to the case $D_{Ch}$.

\satz{Lemma} \label{XF}
Let $X \in Z^p_{eff}(\Pe^t_\C)$, and
$\Pe(W) \subset \Pe^t_\C$ a subspace of codimension $q > t-p$ that
does not meet the support of $X$. Finally $\Pe(F)$ a subspace of
codimension $r \leq t-p$ containing $\Pe(W)$, and intersecting $X$ properly.
Then, for certain constants $c_3,c_6$ depending only on $p,q,r,t$.
\begin{enumerate}

\item
\[ D(\Pe(W),X_F) = D^{\Pe(F)} (\Pe(W), X_F . \Pe(F)) + c_3 \deg X. \]

\item
\[ D^{\Pe(F)}(\Pe(W),X_F . \Pe(F)) = D(\Pe(W),X_F . \Pe(F)) + c_6 \deg X. \]
\end{enumerate}
\end{Satz}

\proof
Assume first $r=t-p$

\vspace{2mm}

1.
In this case the intersection $\Pe(F) . X_F = \Pe(F) . X$ 
is zero dimensional, hence
$X_F$ consists of $t-p$-dimensional subspaces. If $\deg X_F =1$, the
intersction of $X_F$ with $\Pe(F)$ is a single point, hence by \cite{App1}, 
fact 4.8, $D^{\Pe(F)}(\Pe(F). X_F, \Pe(W)) = \log |\Pe(F).X_F,\Pe(W)| + c$ 
where $c$ is a constant only depending on $q,p$ and $t$. Similarly, by the same
fact, $D(\Pe(W),X_F) = \log |\Pe(W),X_F| + c_1$ with $c_1$ only depending
on $r,p$, and $t$. Further, since $X_F$ is orthogonal to $\Pe(F)$ the equality
$|X_F,\Pe(W)| = |\Pe(F).X_F,\Pe(W)|$ holds, and we get
\[ D^{\Pe(F)}(\Pe(W),\Pe(F).X_F) = D(\Pe(W),X_F) + c_3. \]
Since the algebraic distance is additive, for $X_F$ of arbitrary degree
the equality
\[  D(\Pe(W),\Pe(F).X_F) = D(\Pe(W),X_F) + c_3 \deg X \]
follows. 
\vspace{2mm}

2.
Again 
\[ D(\Pe(W),X_F . \Pe(F)) = \log |\Pe(W),X_F . \Pe(F)| + c_4 \deg X, \]
and
\[ D^{\Pe(F)}(\Pe(W), X_F. \Pe(F)) = \log |\Pe(W), X_F. \Pe(F)|+c5 \deg X, \]
for the same reasons as in 1.
The equality 
\[ D^{\Pe(F)}(\Pe(W),X_F . \Pe(F)) = D(\Pe(W),X_F . \Pe(F)) + c_6 \deg X \]
follows.

\vspace{2mm}

Let now $r \leq t-p$.

1.
Let $\Pe(V) \subset \Pe(F)$ be a subspace of codimension $r=t-p$ that contains
$\Pe(W)$, and intersects $X$ properly. 
Since, $(X_F)_V = X_V$, by Lemma \ref{raumraum}, 
\begin{equation} \label{XFa}
D(\Pe(W),X_F) = D(\Pe(W),X_V) + D(\Pe(V),X_F) - D(\Pe(V),X_V).
\end{equation}
On the other hand, consider $(X_F. \Pe(F))_V$ inside $\Pe(F)$. Again, by Lemma
\ref{raumraum}
\[ D^{\Pe(F)}(\Pe(W),X_F. \Pe(F)) = \]
\begin{equation} \label{XFb}
D^{\Pe(F)}(\Pe(W),(X_F . \Pe(F))_V) + D^{\Pe(F)}(\Pe(V),(X_F . \Pe(F))) - 
D^{\Pe(F)}(\Pe(V),(X_F. \Pe(F))_V). 
\end{equation}
We want to show that the left hand side of (\ref{XFa}) equals a constant 
times $\deg Z$ plus the left hand side of (\ref{XFb}). 
Since, $(X_F, \Pe(F))_V = X_V . \Pe(F)$, the first terms on 
the right hand sides coincide modulo a constant times $\deg Z$ by the Lemma
for $r=t-p$. The third terms on the right hand sides are constants times
$\deg Z$ by the proof of \cite{BGS}, Proposition 5.1.1. Finally the 
second terms on the right hand sides coincide modulo a constant times 
$\deg Z$ by \cite{App1}, Lemma 4.13.1.

\vspace{2mm}

2. 
With $\Pe(V)$ as in part one, we have by Lemma \ref{raumraum}
\[ D(\Pe(W),X_F. \Pe(F)) = \]
\begin{equation} \label{XFc}
D(\Pe(W),(X_F . \Pe(F))_V) + D(\Pe(V),(X_F . \Pe(F))) - 
D(\Pe(V),(X_F. \Pe(F))_V). 
\end{equation}
The terms on the right hand side of (\ref{XFc}), and (\ref{XFa}) can be
compared completely analogously as in part one, again using the Lemma
for $r=t-p$, the proof of \cite{BGS}, Proposition 5.1.1, and
\cite{App1}, Lemma 4.13.1.

\proof {\sc of Proposition \ref{raumraum}:}
Since $X_F . \Pe(F) = X . \Pe(F)$, the Proposition simply follows from the two
Lemmata.

\section{Affine Differentiation}

Let $Z$ be a K\"ahler manifold of dimension $d$, i.\@ e.\@ a complex manifold 
equipped with a metric on the tangent space $T_z Z$ for every
$z \in Z$ such that the fundamental form defined by this metric is closed.



A map smooth map $\R^m \to \R^n$ is called analytic, if its Taylor series
locally converges. Clearly if $f:\C^m \to \C^n$ is holomorphic, the induced
map $f_\R: \R^{2m} \to \R^{2n}$ is analytic, and if $\C^n \to \C$ is 
holomorphic the maps $F = |f|: \R^{2n} \to \R$ and $\log F$ are analytic.

Let now $U_\theta$ be a neighbourhood of $\theta \in Z$ as above, and
\[ \varphi, \psi: U \to U_\theta \]
holomorphic charts centerd at the origin.  
Further, denote by
$|\cdot,\cdot|$ the distance on $Z$ as well as the standard distance on
$\A^d(\C) = \C^d$. If $U_\theta$ is relatively compact, then there are
positive constants $c_1,c_2$ such that for every $z_1,z_2 \in U_\theta$,
\begin{equation} \label{abstabsch}
|\varphi^{-1} z_1,\varphi^{-1} z_2| \leq c_1 |z_1,z_2|, \quad
|z_1,z_2| \leq c_2 |\varphi^{-1} z_1,\varphi^{-1} z_2|,
\end{equation}
and the same with $\psi$.

\satz{Lemma}
Let $\partial$ denote the vector
$(\partial/\partial x_1, \partial/\partial y_1, \ldots,
\partial/\partial x_t, \partial/\partial y_t)$, which will also be denoted
$(\partial_1, \ldots, \partial_{2d})$ shortly, and by
$I = (i_1, \ldots, i_{2d})$ a multiindex of order $S$, i.\@ e.\@
$|I|=i_1+ \cdots + i_{2d} = S$. Then, with $\partial^I$ the derivative
$\partial^{i_1}/(\partial x_1)^{i_1} \cdots \partial^{i_{2d}}/
(\partial y_d)^{i_{2d}})$,
\[ \log |\partial^I (\varphi^*f)(0)| \leq 
   \sup_{s \leq S, |J| =s} \log |\partial^J (\psi^* f)(0)|
   + c S \log S, \]
and
\[ \log |\partial^I (\psi^*f)(0)| \leq 
   \sup_{s \leq S, |J|=s} \log |\partial^J (\varphi^* f)(0)|
   + c S \log S, \]
with $c$ a constant depending only on $d$ and the charts $\varphi$, and $\psi$.
\end{Satz}

\proof
Since $(\varphi^{-1} \circ \psi)^*: \C^d \to \C^d$ is holomorphic, the 
induced map $(\varphi^{-1} \circ \psi)_\R^*: \R^{2d} \to \R^{2d}$ is
analytic, hence
\[ \log |\partial^J (\varphi^{-1} \circ \psi)_\R^*(0)| \leq C s \log s, \]
with $C$ only depending on $\varphi$, and $\psi$. Successively applying the
chaine rule gives the desired result.

Let now $M$ be a set of functions $Z \to \R$ closed under sums and 
differences. A grading on $M$ is a map $\deg: M \to \Z$ such that
$\deg (f+g) = \deg f + \deg g$ and $\deg(-f) = - \deg f$ for every
$f,g \in M$.
 
\satz{Definition}
A graded set $M$ of 
smooth functions $f: U_\theta \to \R^{>0}, i \in I$ is said to have a 
holomorphic model if there is an analytic function $g: U_\theta \to \R^{>0}$, 
and for every $f \in M$ a holomorphic function $F: U_\theta \to \C$ such that
\[ f = \log |F| + \deg f \; \log g. \]
\end{Satz}

For multiindizes $I = (i_1,\ldots,i_{2d}),J=(j_1, \ldots, j_{2d})$
write $J \prec I$ iff $j_k \leq i_k$ for all $k=1, \ldots, 2d$. If $J \prec I$
the multiindex $I-J$ is defined.
Further, for $I=(i_1, \ldots, i_{2d}$ a multiindex, and 
$\partial^I$ the corresponding differential operator, write $\partial^I_z$
for the differential operator \\
$\partial^{i_1}/(\partial z_1)^{i_1} \partial^{i_2}/(\partial z_1)^{i_2}
\partial^{i_3}/(\partial z_2)^{i_3} \cdots 
\partial^{i_{2d-1}}/(\partial z_d)^{i_{2d-2}}
\partial^{i_{2d}}/(\partial z_d)^{i_{2s}}$.

\satz{Lemma} \label{CR}
If $\varphi: U \to U_\theta$ is a holomorphic chart, and $M$ a set
of functions that has a holomorphic model, then for every $f \in M$
with holomorphic model $F$, the function $F$ 
locally has a square root $h$, and
\[ \sup_{|I| \leq S}|(\partial^I (\varphi^*f)(0)| \leq
   \sup_{|J| \leq S} \left| \left( \partial^J_z \varphi^* h\right) (0) 
   \right|^2 + O((S+ \deg f) \log (S \deg f)), \]
where the constants implied by the $O$-notation depend on the choic of
the holomorphic charts only.
\end{Satz}

\proof
For $s \leq S$, and $I$ a multiindex of order $s$, 
\begin{equation} \label{CR1}
(\partial^I \varphi^*f)(0) = \sum_{J \prec I}
(\partial^J \varphi^* |F|)(0) \partial^{I-J} \varphi^*g^{deg f})(0). 
\end{equation}
Since $g$ is an analytic function,
\[ |\partial^{I-J} \varphi^*g^{deg f})(0)| \leq 
   (|I-J| \deg f)^{c (|I-J| + \deg f)}  \leq (s \deg f)^{c (s+ \deg f)} \]
for some constant $c$. Hence, the absolute value of (\ref{CR1}) 
is less or equal
\begin{equation} \label{CR2}
(2d)^s \sup_{J \prec I} |(\partial^J \varphi^* |F|)(0)|
(s \deg f)^{c (s+ \deg f)}. 
\end{equation}
Next, with $H = \varphi^* h$,
\[ |(\partial^J \varphi^* |F|)(0)| = |(\partial^J (H\bar{H})(0)|, \]
and the Cauchy-Riemman differential equations imply
\[ \frac{\partial H}{\partial z_i} = \frac{\partial H}{\partial x_i} = 
   -i \frac{\partial H}{\partial y_i}, \]
consequently,
\[ \left|\partial^J H \right| = 
   \left|(-i)^{i_2+ i_4 + \cdots i_{2d}} \partial^J_z H \right| =
    |\partial^J_z H|, \]
and similarly
\[ |\partial^J \bar{H} | = 
   | \overline{\partial_z^J H}|. \]
This implies 
\[ \left|\partial^J |\varphi^* f| \right| =
   \left|\partial^J (H \bar{H}) \right| \leq
   2^{|J|} \sup_{K \prec J} \partial_z^K H  \overline{\partial_z^{J-K} H} \leq
   2^s \sup_{|K| \leq s} |\partial_z^K H|^2. \]
   
Hence, (\ref{CR2}) is less or equal
\[ (4d)^s (s \deg f)^{c(s+ \deg f)} \sup_{|K| \leq s} |(\partial_z^K)(0)|^2, \]
proving the Lemma.

\vspace{2mm}

By standard complex ananlysis, the derivatives of a function $f$
at $\theta$ that has a holomorphic
model $F$ can be estimated by the values of $F$ on $U_\theta$ in 
the following way.

\satz{Proposition} \label{Cauchy}
Let $Z$ be a K\"ahler manifold of dimension $d$, and $\theta \in Z$
a point with 
neighbourhood $U_\theta$; further $\varphi: U \to U_\theta$ an affine chart, 
and $f$ a smooth function on $U_\theta$ that has a holomorphic
model $F$. Then, for every $S \in \N$, every $I$ with $|I|=S$, and
every $R \in \R$ such that the ball of radius $R$ around $\theta$ is contained
in $U_\theta$,
\[ \log |(\partial^I \varphi^*f) (0)| \leq 
   \sup_{|\theta,z| \leq R} \log |F(z)| - 2S \log R + O(S \log S). \]
\end{Satz}

\proof
Let $h$ be a local square root of $F$ and $H = \varphi^* h$.
Lemma \ref{CR} implies for any multiindex with $|I|=S$,
\[ \log |(\partial^I \varphi^* f)(0)| \leq
   \sup_{|J|\leq S}|(\partial^J_z H)(0)|^2 = \]
\begin{equation} \label{Cau4} 
\sup_{|J| \leq S} \log
   \left| \left( \frac{\partial^s}{(\partial z_1)^{j_1+j_2} \cdots
         (\partial z_d)^{j_{2d-1}+ j_{2d}}}H \right)(0) \right|^2+
   (S + \deg f) \log (S \deg f). 
\end{equation}
By the multidimensional Cauchy formula,
\[
\left| \left( \frac{\partial^s}{(\partial z_1)^{j_1+j_2} \cdots
         (\partial z_d)^{j_{2d-1}+ j_{2d}}}H \right)(0) \right|^2 \leq \]
\[ \frac1{(2\pi)^d} \left| \int_{z_1=R'} \cdots \int_{z_d = R'}
\frac{H}{z_1^{k_1+k_2} \cdots z_d^{k_{2d-1}+ k_{2d}}} \; 
d z_1 \cdots d z_d \right|^2 \leq 
\frac{1}{2\pi^d} (R')^{-2s} \sup_{|z,0| \leq R'} |H(z)|^2 \]
with $R' = R/c_2$, and $c_2$ from (\ref{abstabsch}),
which in turn equals
\[ \frac{1}{2 \pi^d}(R')^{-2s} \sup_{|z,\theta|\leq R} |F(z)|. \]
Inserting this into (\ref{Cau4}) finishes the proof.

\subsection{Projective space}





\satz{Lemma}  \label{proj}
Let $\theta \in \Pe^t(\C)$.

\begin{enumerate} 


\item
If $f \in \Gamma(\Pe^t,O(D))$ is a global section whose restriction
to $\theta$ is nonzero, then the function
\[ \zeta \mapsto \log |f_\zeta|, \]
hence likewise the function
\[ D(\zeta,\di f) = \log |f_\zeta| - \int_{\Pe^t} \log |f| \; \mu^t \]
have holomorphic models.

\item
For $U_\theta$ the circle of radius $r<1$ around $\theta$,
the inequalities (\ref{abstabsch}) hold with $c_1=1/\sqrt{1-r^2}, c_2=1$.

\end{enumerate}
\end{Satz}

\proof
1.
If the global sections 
$\Gamma(\Pe^t,O(D))$ of the line bundle $O(D)$ are given by homogeneous
polynomials of degree $D$ in variables $z_0, \ldots, z_t$ with
$z_1(\theta)= \cdots = z_t(\theta) =0$, the map
$f(z_0,\ldots,z_t) \mapsto F(z_1, \ldots,z_t)=f(1,z_1,\ldots,z_t)$ maps 
$\Gamma(\Pe^t,O(D))$ to the space of polynomials on $\A^t$ of degree
at most $D$. Further, with $\varphi: \A^t \to \Pe^t$ the affine chart
with $\varphi(0) = \theta$ and
$\zeta = (\zeta_0, \ldots, \zeta_t)$ in $\varphi(\A^t(\C))$,
\[ |f_{\varphi \zeta}| = \frac{|F(\zeta_1, \ldots, \zeta_t)|}
   {|(1,\zeta_1, \ldots, \zeta_t)|^D}. \] 
Since $g:(z_1,\ldots,z_t) \to |(1,z_1, \ldots, z_t)|$ is analytic, it follows 
that $F(\zeta_1, \ldots, \zeta_t)$ is a holomorphic model for 
$\zeta \to \log |f_\zeta|$, and 
$F(\zeta_1, \ldots,\zeta_t) 
\exp \left( - \int_{\Pe^t} \log |f| \; \mu^t \right)$
is a holomorphic model for $D(\zeta,\di f)$.

\vspace{2mm}

2.
With $\zeta = \varphi(z)$ we have $|0,z| = |z|$, and
\[ |\theta,\zeta| = \sqrt{\frac{|(z_1, \ldots, z_t)|^2}
    {1+{|(z_1, \ldots, z_t)|^2}}}, \]
implying the claim.



\subsection{Grassmannians}

Let now $G(\C) = G_{p,t+1}(\C)$ be the Grassmannian of $p$-dimensional 
subspaces of $\C^{t+1}$. On $G$, there is the line bundle $L$ defined as the
determinant of the canonical quotiont bundle. Further, there is the canonical
harmonic $(1,1)$-form $\mu_G = c_1(\bar{L})$.
This metric explicetely can be described as follows:

Let $W,W' \in G(\C)$, and $S_{W'}$ the unit spere in $W'$ Then, 
\[ |W,W'| = \sup_{w \in S_{W'}} |pr_{W^\bot}(w)|, \]
where $pr_{W^\bot}$ is the orthogonal projection to the orthogonal complement
of $W$.

Let $W_0$ be any $p$-dimensional subspace of $\C^{t+1}$. There is the 
following holomorphic chart $\varphi:\A^{(t+1-p)p} \to G$:
Let $w_1, \ldots, w_p$ be an orthonormal basis of $W_0$, $V= W_0^\bot$,
and $U^-$ the unipotent radical of the subgroup of $GL(\C^{t+1}$ that
leaves $V$ invariant. Then the big cell in the Bruhat decomposition of
$G_{p,t+1}$ centered at $W_0$ consists of the subspaces
$u W_0, u \in U^-$. The map
\[ \varphi: \A^{t+1-p)p} \cong U^- \to G, \quad
   u \mapsto u W_0 \]
is certainly holomorphic.

\satz{Lemma} \label{gras}
 
\begin{enumerate}


\item
For any hypersurface $Z = \di f$ in $G$ the map
\[ G \to \R, \quad V \mapsto
   D(V, Z) = \log |f_V| - \int_G \log |f| \mu_G^{p(t+1^-p)} \]
has a holomorphic model.

\item 
The inequality (\ref{abstabsch}) holds with $c_2=1$, and $c_1$ some constant
depending on $p$ and $t$.

\end{enumerate}
\end{Satz}


\proof
1.
Let $\check{E}$ be the vectorbundle on $G_{p,t+1}$ that attaches to
each point $W \in G$ the dual vector space $\check{W}$ of $W$.
The global sections of $\check{E}$ are the vectors 
$\check{v} \in (\C^{t+1})\check{}$,
Since, $\CL_G = \Lambda^p \check{E}$ the global sections of
$\CL_G^{\otimes D}$ are symmetric produckts of vectors of the form
$\check{v}_1 \wedge \cdots \wedge \check{v}_p, \check{v}_i \in 
(\C^{t+1})\check{}$. Now if
\[ f = \prod_{j=1}^D \check{v}_{1j} \wedge \cdots \wedge \check{v}_{pj} \]
is such a global section, then for $\tilde{W}$ not in the support of $\di f$,
\[ D(\di f, \tilde{W}) = \log |f_{\tilde{W}}| - 
   \int_G \log |f| \mu_G^{p(t+1)}. \]
Further if $w_1, \ldots, w_p$ is an orthonormal Basis of $W_0$, and $U^-$
is as defined above, then for $W = u W$,
\[ |f_{\tilde{W}}| = 
   \frac{|\det ((\prod_{j=1}^D \check{v}_{ij})
   (u w_i)^{\otimes D})_{i=1, \ldots, p,k=1, \ldots,p}}
                {|u w_1 \wedge \cdots \wedge u w_p|^D}, \]
and the function inside the numerator is certainly a holomorphic function
of $u \in U^- \cong \A^{p(t+1)}$.

\vspace{2mm}

2. Is clear.


\vspace{2mm}





\satz{Lemma} \label{urhol}
Let $p+q \leq t+1$, and $\Pe(W)$ a subspace of dimension $p-1$, and
$U_W$ a  neighbourhood of $W$ in $G_{p,t+1}$. Then, if
$LZ^{t+1-q}(\Pe^t) \subset Z^{t+1-q}_{eff}(\Pe^t)$ is the subgroup generated 
by subspaces of dimension $p-1$ that intersect each $\Pe(\tilde{W})$ with
$\tilde{W} \in U_W$ properly, the set of functions
\[ f_{Z}: U_W \to \R, \quad \tilde{W} \mapsto D(Z,\Pe(\tilde{W})), \quad
   Z \in LZ^{t+1-q}(\Pe^t) \]
has a holomorphic model.
\end{Satz}

\proof
For $\Pe(V) \subset \Pe^t$ a subspace of dimension $p-1$ intersecting
each $\Pe(\tilde{W}$ with $\tilde{W} \in U_W$ properly,
let $\check{F}_V$ be the vector bundle on $U_W$ that attaches to
each $\tilde{W}$ the space $(\tilde{W}/(V \cap \tilde{W}))^{\check{}}$,
and the line bundle $\CL_V = \bigwedge^{t+1-q}$.
Further, let 
\[ \check{v}_1, \ldots, \check{v}_{t+1-q} \in (\C^{t+1})^{\check{}} \]
be linear forms that are orthonormal and zero on $V$, and define
$f^V = 
\check{v}_1 \wedge \cdots \wedge \check{v}_{t+1-q} \in \Gamma(U_W, \CL)$.
If $P(V)$ is the parabolic subgroup of $Gl(\C^{t+1})$ that leaves 
$V$ invariant, then the group $U^- \cap P(V)$ operates transitively
on $U_W$, and for $\tilde{W} = u W$ we have
\[ |f^V_{\tilde{W}}| =
   \frac{|\det (\check{v}_i (u w_k))_{i=1, \ldots, t+1-q, k=1, \ldots, t+1-q}|}
        {|u w_1 \wedge \cdots \wedge w_{t+1-q}|}. \]
The formula inside the absolute value in the numerator is linar from
in $u$. Furhter, the above expression
equals the sine of the angle between $V$ and $\tilde{W}$, and 
by \cite{App1}, Fact 4.8 $D(\Pe(V), \Pe(\tilde{W}))$, modula a fixed
constant, eqals
the logarithm of the sine of the angle between $V$ and $\tilde{W}$.
equals. Now, for $Z \in LZ^{t+1-q}$ arbitrary,
$Z = \sum_{V} n_V V$, and the function
\[ \prod_V (f^V)^{n_V} \]
models the function $D(Z,\Pe(\tilde{W}))$.

\satz{Lemma} \label{tub}
Let $X \subset \Pe^t_\C$ be a subvariety of dimension $p$, and
$G = G_{t,t-p}$.
\begin{enumerate}

\item
There is a positive constant $c_1$, only depending on $t$, and $p$ such that
\[ \mu(V_X) = \int_{V_X} \mu_G^{(t+1-p)p-1}  = c_1 \deg X. \]

\item
Let $U_\epsilon(V_X)$ be the tubular neighbourhood
\[ U_\epsilon(V_X) := \{ V \in G(\C) | |V,V_X| \leq \epsilon \} \]
of $V_X$. Then, there is a positive constant, depending only on $t$, and $p$
such that
\[ \mu_G(U_\epsilon(V_X)) = \int_{U_\epsilon(V_X)} \mu_G^{(t+1-p)p} \leq
   c_2 \deg X \epsilon. \]
\end{enumerate}
\end{Satz}

\proof
1.
Since $\mu_G$, and thereby $\mu_G{(t+1-p)p-1}$ are closed forms, 
the integral 
\[ \int_{V_X} \mu_G^{(t+1-p)p-1} \]
depends only on the cohomology class of the cycle $V_X$, and thereby only
on the class of $V_X$ in $CH^1(G)$. Since in this last group
\[ [V_X] = \deg X [V_{\Pe(W)}], \]
where $\Pe(W) \subset \Pe^t$ is any projective subspace of codimension
$p$, with \\ $c_1:=\int_{V_{\Pe(W)}} \mu_G^{(t+1-p)p-1}$ the equality
$\int_{V_X} \mu_G^{(t+1-p)p-1}  = c_1 \deg X$ follows.

\vspace{2mm}

2.
Since $G$ has positive curvature, this immediately follows from part one.

\satz{Lemma} \label{trfunk0}
Let $\Pe(W) \subset \Pe(F)$ be subspaces of $\Pe^t$ of codimensions
$q \geq r$. Every subspace $\Pe(F')$ of codimension $r$ contains a subspace
$\Pe(W')$ of codimension $q$ such that
\[ |W,W'| \leq |F,F'| \]
as points in the corresponding Grassmannians.
\end{Satz}

\proof 
The proof is elementary linear algebra. Let $pr_F: \C^{t+1} \to F$ be
the orthogonal projection to $F$, and define $W'$ as the intersection
$pr_F^{-1}(W) \cap F'$. Then $W'$ has codimension $q$, is contained in $F'$,
and every vector $w \in W'$ may be written as $W= w_1 + v_1$ with
$w_1 \in W, v_1 \in F^\bot \subset W^\bot$. Hence,
\[ pr_{W^\bot}(w) = v_1 = pr_{F^\bot}(w), \]
and  since $W' \subset F'$,
\[ |W,W'|  = \sup_{w \in W'} |pr_{W^\bot}(w)| =
   \sup_{w \in W'} |pr_{F^\bot}(w)| \leq \sup_{v \in F'} |pr_{F^\bot}(v)| =
   |F,F'|. \]

\vspace{2mm}

Next, there is the following functoriality for differential operators
on Grassmannians of subspaces of different dimension.
Let $p+q \leq t+1$, $\Pe(W) \subset \Pe^t(\C)$ be a subspace of dimension 
$q-1$, and $U_W$ an open subset of $W$ in $G_{q,t}$ that is contained
in the big cell in the Bruhat decomposition centered at $W$.

For $\Pe(F)$ a subspace of dimension $p+q-1$ containing $\Pe(W)$, 
let $V$ be the orthogonal complement of $W$ in $F$, and define the map
\[ f: U_W \to G_{p+q,t}(\C), \quad
   \tilde{W} \mapsto \tilde{W} \oplus V. \]
Clearly, $f$ is a holomorphic map.


\satz{Lemma} \label{trfunk}
Let $| \cdot , \cdot |$ be the canonical metric on the Grassmannian. With
the above notations, 

\begin{enumerate}

\item
\[ \forall \; \tilde{W} \in U_W \;
   |F,\varphi(\tilde{W})| \leq |W,\tilde{W}|. \]

\item
Let $\varphi: \A^{} \to U_w, \psi:\A^{} \to G_{}$. Then,
\[ \sup_{|I|=s} |\partial_z^I (\psi^{-1} \circ f \circ \varphi)(0)| \leq
   c s \log s, \]
where $c$ is a constant depending only on $p,q$, and $t$ but not on the 
choice of $W$ and $F$.




\end{enumerate}
\end{Satz}

\proof
1. Let $\tilde{W} \in U_W$ and $\tilde{F} = \varphi(\tilde{W})$.
Since $V$ is contained in $F$ , we have
\[ |F,\tilde{F}| =  \sup_{v \in S_{\tilde{F}}} |pr_{F^\bot}(v)| =
   \sup_{w \in S_{\tilde{F} \cap V^\bot}} |pr_{F^\bot}(v)|. \]
Now let $v \in S_{\tilde{F} \cap V^\bot}$ be a vector where value of the last
supremum is achieved, and $w \in S_{\tilde{W}}$ be a vector such that
$|W, \tilde{W}| = |pr_{W^\bot}(w)|$. We have to show
\[ |pr_{F^\bot}(v)| \leq |pr_{W^\bot}(w)|, \]
which again boils down to an elementary calculation in linear algebra.

\vspace{2mm}

2. This is an immediate consequence of the holomorphicity of $f$, and the
fact that $\psi^{-1} \circ f \circ \varphi$ is the same for all $W,F$ modulo
a transformation by a $g \in SU(t+1)$.



\section{The derivated algebraic distance}

\subsection{Hypersurfaces and points}

Let $Z$ be a regular projective algebraic variety of dimension $d$
over $\C$, and fix a 
K\"ahler structure on $Z$. For $f$ a global section of some line bundle on 
$Z$, $X = \di f$ an effective cycle of pure codimension $1$ on $Z$, and
$\theta \in Z$ a point not contained in the support of $X$ the algebraic
distance $D(\theta,X)$ equals
\[ D(\theta,X) = \log |f_\theta| - \int_Z \log |f| \mu^d =
   - \frac12 \int_X g_\theta + \frac12 \int_Z g_\theta \mu, \]
where $\mu$ is the chosen K\"ahler form on $Z$, and $g_\theta$ is a
green form of log type for $\theta$. (See \cite{App1}, Definition 4.1).

For $Y$ a cycle in $Z$ of pure codimension $1$ the function
$D(Y,\theta)$ clearly is smooth in a neighbourhood of $\theta$. 
This leads to the following definition.

\satz{Definition}
With the above notations, let $Y \in Z_{eff}(Z)$ be an effective cycle
of pure codimension $1$ in $Z$ on
$Z$, and $\theta \in Z$ a point not contained in the support of $Z$.
Further $\varphi: U \subset \A^d \to U_\theta$ an affine chart.
The derivated algebraic distance of $Z$ to $\theta$ is defined as
\[ D^S (\theta, Z) := \sup_{|I| \leq S} 
   \log  \left|\left(\partial^I \exp (D(z,Z))\right)(\theta) \right|. \]
\end{Satz}

\satz{Proposition} \label{grashol}
Let $q \leq t$, and $G=G_{q,t}$ the Grassmannian of $q$-dimensional subspaces
of $\C^{t+1}$. Then, for effective cycle  $Z$ of codimension $q$
in $G$and every point 
$W \in G$ not contained in $Z$, the algebraic distance $D(Z,W)$ has a 
holomorphic model.
\end{Satz}

\proof
Since $CH^1(G) \cong \Z$, there is some global section $f$ of the
canonical line bundle $L_G$ on $G$ such that $Z = \di f$.
The Proposition thus follows form Lemma \ref{gras}.1.



\satz{Corollary}
With $f$ a global section of a line bundle $L$ on $Z$, $Y = \di f$, 
$\theta$ a point not contained in the support of $Y$, 
$\varphi: U \to U_\theta$ an affine chart, and $H$ a local square root
of $\varphi^* f$,
\[ D^S(Y,\theta) = 
   \sup_{|J| \leq S} \left| \left( \frac{\partial^s}
   {(\partial z_1)^{j_1} \cdots (\partial z_t)^{j_t}} H \right)
   (0) \right|^2 + O(S \log S S). \]
\end{Satz}

\proof
Follows from the Proposition and Lemma \ref{CR}.

\vspace{2mm}

The derivated algebraic diestance of a hypersurface in $\Pe^t$ can also be
estimated against the values of the derivations of the global sections
directily as stated in Theorem \ref{hypeb}.

\proof {\sc of Theorem \ref{hypeb}}
By the proof of Lemma \ref{proj},
with $\varphi: \A^t \to \Pe^t$ a homogeneous chart centered at $\theta$,
and $\zeta_1, \ldots, \zeta_t$ the coorindates in $\A^t$
\[ D(\di f, \varphi(\zeta)) = \log |F(\zeta)| - \int_{\Pe^t} |f| \mu^t - 
   D \log |(1,\zeta_1,\ldots,\zeta_t)|, \]
hence
\[ D^S(\di f, \theta) = \sup_{|I| \leq S} \log 
   |(\partial^I |F|)(0)| - \int_{\Pe^t} |f| \mu^t + O(S \log S). \]
By the first equality above the Proposition holds for $S=0$. Assume
that for a multiindex $I$ with $|I|=S-1$,
\[ |(\partial^I |F|)(0)| \leq |(\partial_z^I F)(0)|. \]
Then, for any $j=1, \ldots, t$,
\[ |(\partial/\partial x_j \partial^I |F|)(0)| \leq
   |(\partial/\partial x_j |(\partial_z^I F)(0)|. \]
With $G = \partial_z^I F$, and $G_r,G_i$ the real and imagainary parts of $G$ 
this equals
\[ \left| \frac{G_r \partial G_r/\partial x_j + G_i \partial G_i/\partial x_j}
   {\sqrt{G_r^2+G_i^2}} (0)\right| \leq
   \left| \frac{\partial G_r}{\partial x_j} + \frac{\partial G_i}{\partial x_j}
   (0)\right| \leq \left|\frac{\partial G_r}{\partial x_j}(0)\right| + 
   \left|\frac{\partial G_i}{\partial x_j}(0)\right|, \]
which by the Cauchy-Riemman-equations equals
\[ 2 \left|\frac{\partial G}{\partial z_j}\right| =
   2 \left|(\partial/\partial z_j \partial^I F)(0)\right|. \]
The Proposition follows by complete induction.

\subsection{Effective cycles in projective space and projective subspaces}

\satz{Proposition} \label{prohol}
Let $p,q \in \N$, $Z \in Z^p_{eff}(\Pe^t_\C)$, and $\Pe(W) \subset \Pe^t(\C)$
a subspace of codimension $q$ that is regular with respect to $Z$. Let further
$G=G_{t+1-q,t}$ denote the Grassmannian of $q$- codimensional subspaces of 
$\Pe^t(\C)$. There is an open neighbourhood $U_W$ of $W$ in $G$ such that 
every $\tilde{W} \in U_W$ is regular with respect to $Z$. Further, 
\begin{enumerate}

\item
If $t \leq p+q \leq t+1$, or more generally if $p+q \leq t+1$, and for
every $\tilde{W}$ in some neighbourhood $U_W$ of $W$ the intersection
$Z . \Pe(\tilde{W})$ is a sum of projective subspaces, then
the function $D(Z,\Pe(\tilde{W})$ has a holomorphic 
model, in particular is a smooth function.

\item
If $p+q > t+1$, the function $D(Z,\bullet)$ is smooth on $U_W$.

\end{enumerate}
\end{Satz}

This leads to the following Definition:

\satz{Definition}
Let $Z \in Z_{eff}(\Pe^t_\C)$ be an effective cycle in projective space, and
$\Pe(W) \subset \Pe^t(\C)$ a projiective subspace of codimension $q$
that is regular with respect to $Z$. 
For $W$ the corresponding point in the Grassmanninan
$G=G_{t+1-q,t+1}$ there is a connected simply connected neigbourhood $U_W$ of 
$W$ in $G$ such that for every $\tilde{W} \in U_W$ the space $\Pe(\tilde{W})$ 
is likewise regular 
with respect to $Z$. Let $\varphi^*U \to U_W$ be the affine chart from 
section 3, and define the derivated algebraic distance of $\Pe(W)$ to $Z$ as
\[ D_{\bullet}^S (\theta, Z) := \sup_{|I| \leq S} 
   \log  \left|(\partial^I \exp \left(\varphi^* D_\bullet(z,Z)\right))(0) 
   \right|, \]
where $D_\bullet(\theta,Z)$ either denotes the algebraic distance 
$D_G(\theta,Z)=D_1(\theta,Z)$ or the algebraic distance $D_{Ch}(\theta, Z)$
of $\theta$ to $Z$ as defined in \cite{App1}, section 4.
\end{Satz}

\subsubsection{Points}

\proof {\sc of Proposition \ref{prohol}.2 for $p=t$:}
For $Z \in Z^t(\Pe^t_\C)$, and $\Pe(W) \subset \Pe^t$ a subspace of
codimension $q \leq t-1$ that does not meet the support of $Z$, let
\[ Z = \sum_{i=1}^{\deg Z} z_i, \]
and $x_i, i=1, j=1, \ldots, q$ global section of
$O(1)$ of length $1$ on $\Pe^t$ such that
\[ \Pe(W) = \di x_1 \cap \ldots \cap \di x_q. \]
For $U_W$ a neighbourhood of $W$ in $G_{t-q,t+1}$, such that for every
$V \in U_W$, no $z_i$ lies in $\Pe(V)$, and define
\[ g_i : U_\theta \to \R, \quad z \mapsto 
   \frac{x_1(z_i) \overline{x_1(z_i)} + \cdots 
   x_q(z_i) \overline{x_q(z_i)}}
   {1+x_1(z_i) \overline{x_1(z_i)} +\cdots x_q(z_i)\overline{x_q(z_i)}}. \]
Clearly $g_i$ is smooth and nonzero on $U_W$. Hence, the
square root $f_i$ of $g_i$ is also smooth on $U_W$.
Since, $|z,z_i| = |g_i(z)|$, and by 
\cite{App1}, Fact 4.8, there is constant $c$ such that
\[ D(Z,\Pe(W)) = c \deg Z + \sum_{i=1}^{\deg Z} \log |z,z_i|, \]
the claim follows.

\vspace{2mm}

For zero dimensional cycles $X$ the derivated
algebraic distance to a ponint $\theta$
not on $X$ takes a particularily simple form

\satz{Proposition} \label{Punkte}
\begin{enumerate}

\item
Let $Z \in Z^t_{eff}(\Pe^t)_{\bar{\C}}$ be an effective cycle of pure dimension
zero, and $\theta \in \Pe^t_\C$ a point not contained in the support of
$Z$. If $Z = \sum_{i=1}^{\deg Z} z_i$ is ordered in such a way that
$|z_1,\theta| \leq \cdots \leq |z_{\deg Z},\theta|$, then for every 
$S \leq \deg Z$
\[ D^S(Z,\theta) \leq \sum_{i=S+1}^{\deg Z} \log |z_i,\theta| + 
   O(S \log \deg Z), \]
and for every $S \leq \deg Z/3$
\[ 2\sum_{i=S+1}^{\deg Z} \log |z_i,\theta| \leq
   D^{3S}(Z,\theta) + O((S+ \deg Z) \log \deg Z). \]

\item
For $p \leq t$ let $Z \in Z^p_{eff}(\Pe^t)$ such that
$Z = \sum_{i=1}^{\deg Z} \Pe(W_i)$ with each $\Pe(W)$ a projective subspace,
and $\theta$ a not contained in $\Pe(W_i)$ for all $i=1, \ldots, \deg Z$. Then,
for every $S \leq \deg Z$,
\[ D^S(\theta,Z) \leq \sum_{i=S+1}^{\deg Z} \log |\theta, \Pe(W_i)| +
   O(S \log \deg Z). \]
\end{enumerate}
\end{Satz}

The next three Lemmata will be proved in the appendix

\satz{Lemma} \label{hilf}
Let $x_1, \ldots, x_n \in [-1,1] \setminus \{0\}$ with
\[ |x_1| \leq |x_2| \leq \cdots \leq |x_n|, \quad \mbox{and}
   \quad f(x) := \prod_{i=1}^n (x-x_i). \]
Then, for $s < n/3$,
\[ \frac1{(2s+1)(3n^3)^{s+1}} \prod_{i=s+1}^n |x_i|^2 \leq
   \sup_{0 \leq j \leq 3s} |f^{(j)}(0)|, \]
and for $s \leq n$,
\[ \sup_{0\leq j \leq s} |f^{(j)}(0)| \leq \frac{n!}{(n-s)!} 
   \prod_{i=s+1}^n |x_i|. \]
\end{Satz}

\satz{Lemma} \label{hilfzwei}
Let $\theta, z_1, \ldots, z_n$ be points in $\Pe^t(\C)$, with
$\theta \neq z_i \forall i=1, \ldots, n$, and
$\varphi: \A^t(\C) \to \Pe^t(\C)$ the affine chart centered at $\theta$.
With
\[ f(z) = \prod_{i=1}^n |z,z_i|, \]
and $F = \varphi^* f$, for every multiindex
$I = (i_1, \ldots, i_t)$ with $|I| = s \leq n/3$, 
\[ \frac1{(2s+1)(3n^2)^{s+1} n^n} \prod_{s+1}^n |\theta,z_i|^2 \leq 
   \sup_{|I| \leq 3s} \left|(\partial^I F)(0) \right|, \]
and for $s \leq n$,
\[ \sup_{|I| \leq s} 
   \left|(\partial^I F)(0)\right| \leq
   n^s \prod_{i=s+1}^n |\theta,z_i|. \]

\end{Satz}

\satz{Lemma} \label{hilfdrei}
Let $\Pe(W_i) \subset \Pe^t_\C, i=1, \ldots n$ be subspaces of fixed 
codimension $p$, and $\theta \in \Pe^t(\C)$ a point contained in none of 
them. Further, $\varphi: \A^t(\C) \to \Pe^t(\C)$ the affine chart centered
at $\theta$ and
\[ f(z) = \prod_{i=1}^n |z, \Pe(W_i)|. \]
Then with $F = \varphi^*f$, and $s \leq n$,
\[ \sup_{|I| \leq s} \left| (\partial^I F)(0) \right| \leq
   n^s \prod_{i=1}^n |\theta, \Pe(W_i)|. \]
\end{Satz}

\proof {\sc of Proposition \ref{Punkte}:}
1. The inequality
\[ D^S(Z,\theta) \leq O(S \log \deg Z) \]
follows immediately from the second inequality of Lemma \ref{hilfzwei}, 
together with the equality
\begin{equation} \label{Punkte1}
D(\theta,x) = \log |x,\theta| + c, 
\end{equation}
with $c$ a constant only depending on $t$ from \cite{App1}, Fact 4.10.
The second inequality follows from the first equality of Lemma \ref{hilfzwei}, 
and again equality (\ref{Punkte1}). 

\vspace{2mm}

2. Follows in the same way using Lemma \ref{hilfdrei}, and the equality
\[ D(\theta, Z) = \sum_{i=1}^{\deg Z} \log |\theta, \Pe(W_i)| +
   c \deg Z, \]
which follows again from \cite{App1}, Fact 4.10 together with the
additivity of the algebraic distance.

\subsubsection{The general case}

\proof {\sc of Proposition \ref{prohol}. 1:}
Assume first that $p+q =t+1$. Then, by the proof of \cite{App1},
Proposition 4.14,
\[ D(Z,\Pe(W)) = D(V_Z, W), \]
with $V_Z$ $G^* F_* Z \in Z^1_{eff}(G_{t+1-p,t})$
from the Correspondence (\ref{grascor})
Thus, the claim in this case is Proposition \ref{grashol}.

\vspace{2mm}

For $p+q \leq t$, and $Z . \Pe(\tilde{W})$ equal to a sum of projective
subspaces for every $\tilde{W}$ in some neighbourood of $W$ assume
first $p+q = t$, and
let $\Pe(V) \subset \Pe^t(\C)$ be a subspace of codimension one 
that does not meet the support $\Pe(\tilde{W}) .Z$ for every 
$\tilde{W} \in U_W$. By \cite{App1}, Proposition 4.12,
\begin{equation} \label{proholB1}
D(Z,\Pe(\tilde{W})) =
D(\Pe(V),\Pe(\tilde{W})) + D(Z, \Pe(V) \cap \Pe(\tilde{W})) -
D(Z . \Pe(\tilde{W}),\Pe(V)). 
\end{equation}
Hence, the claim of the Proposition holds for $D(Z,\Pe(\tilde{W})$ if 
it holds for every term on the left hand side of (\ref{proholB1}).
For the second term, the claim follows from Case A. For the first term,
it follows from Lemma \ref{urhol}.

For the third term, let $U^\circ_W \subset U_W$ be the subset such
that for every $\tilde{W} \in U^\circ_W$ the space $\Pe(\tilde{W})$
does not meet the singular locus of $Z$, and $Z . \Pe(\tilde{W})$ has
no double points. Clearly $U_W \setminus U^\circ_W$ is contained in
the intersection of $U_W$ with the zero set of a holomorphic function. 
Further, let
$X$ be the global section of $O(1)$ such that $\Pe(V) = \di X$. Then,
\[ D(Z . \Pe(\tilde{W}), \Pe(V)) =
   \log \prod_{i=1}^{\deg Z} 
   \frac{|X(z_i)|}{\sqrt{1+ X(z_i) \overline{X(z_i)}}} + c \deg Z. \]
On $U_W^\circ$ the coordinates of $z_1, \ldots, z_{\deg Z}$ are holomorphic
functions of $\tilde{W}$. Hence $X(z_i)$ is holomorphic on $U_W^\circ$.
Since $X(z_i)$ is continious on $U_W$, and $U_W \setminus U^\circ_W$ is 
contained in the intersection of $U_W$ with the zero set of a holomorphic 
function, it follows that $X(z_i)$ is holomorphic on all of $U_W$, finishing
the proof.

If $p+q < t$, let $\Pe(V) \subset \Pe^t(\C)$ be a subspace of codimension
$t+1-p-q$ that does not meet the support of $\Pe(\tilde{W}). Z$ for every
$\tilde{W} \in U_W$. Then (\ref{proholB1}) still holds and the first
and second term on the right hand side by the same arguments have
holomorphic models. Further, by assumption
$Z. \tilde{W} = \sum_{i=1}^{\deg Z} \Pe(\tilde{W})$, hence
by \cite{App1}, p. 28,
\[ D(Z. \Pe(\tilde{W},\Pe(V)) =
   \sum_{i=1}^n D(\Pe(\tilde{W}_i), \Pe(V)) =
   \sum_{i=1}^n D(\tilde{W}_i,V_{\Pe(V)}), \]
where the $\tilde{W}_i$ are points and $V_{\Pe(V)}$ is a hypersurface
in $G_{t-p-q,t+1}$, and
algeraic distance of effective cycles in $\Pe(F)$, 
\[ D(\Z.\Pe(\tilde{W}), \Pe(V)) =
   \log \prod_{i=1} \frac{X(\tilde{W}_i)}
   {\sqrt{1+X(\tilde{W}_i)\overline{X(\tilde{W}_i)}}}, \]
with $X$ a global section of the canonical line bundle of the Grassmannian
such that $V_{\Pe(V)} = \di X$.
Taking again $U_W^\circ \subset U_W$ as the subset sucht that for every
$\tilde{W} \in U_W^\circ$ the space does not meet the singular locus
of $Z$, and $Z. \Pe(\tilde{W})$ has no double points,
the coordinates of each $\tilde{W}_i$ depend holomorphically on 
$\Pe(\tilde{W})$, and one can repeat the argument above.

\proof {\sc of Proposition \ref{prohol} for $p+q >t+1$:}
Let $\Pe(F)$ be a subspace of codimension $t-p$ containing $\Pe(W)$,
and intersecting $Z$ properly. By Proposition \ref{raumraum},
\[ D(Z,\Pe(W)) - D(Z_F,\Pe(W)) = D(Z,\Pe(F))- D(Z_F,\Pe(F)). \]
By part one of the Proposition, the left hand side is smooth.
Further, since \\ $D(Z_F),\Pe(W)) = D(Z_F . \Pe(F),\Pe(W)) + c \deg X$,
and $Z_F . \Pe(F)$ consitst of points, the second term on the left hand
side is smooth by part two of the proposition for $p=t$ proved above.
Hence, $D(Z,\Pe(W))$ is likewise smooth.


\section{Reduction of the derivated algebraic distance to 
derivated algebraic distances to points}

Let $Z \subset \Pe^t_\C$ be an algebraic subvariety of codimension $p$, and
$\theta \in \Pe^t(\C)$ a point not contained in $Z$. 
\cite{App1}, Proposition 4.16 implies the existence of a projective subspace
$\Pe(F) \subset \Pe^t$ of codimension $t-p$ that intersects $Z$ properly
and contains $\theta$, such that
\[ c_1 \deg Z \geq D(\Pe(F),Z) \geq -c_2 \deg Z, \]
with positive constants $c_1,c_2$ only depending on $p$, and $t$, and, 
if $D^{\Pe(F)}(\bullet,\bullet)$ denotes the
\[ D(Z,\theta) = D^{\Pe(F)}(Z. \Pe(F), \theta) + O(\deg Z). \]
The existence of such a space thus allowed to reduch teh algebraic distance of
a point $\theta$ to an effective cycle $Z$ to the algebraic cycle to the 
algebraic distance of $\theta$ to $Z.\Pe(F)$ which is zero dimensional.
For decomposing the derivated algebraic distance however, the condition 
$D(Z, \Pe(F)) \geq \ - c_2 \deg Z$ is not good enough, because derivatives
of a function may be very small or big even if the values of the function
are not; to assure that there is a space $\Pe(F)$ such that the
derivations of $\exp D(\Pe(F),Z)$ and $\exp (-D(\Pe(F),Z))$ are also
bounded in terms of $\det Z$, one has to look for a space that contains a
smaller subspace that has not too small distance to $Z$.

\satz{Theorem} \label{zerl}
Let $p \leq t$, and $Z \subset \Pe^t(\C)$ be an effective cycle of pure 
codimension $p$; further $\theta \in \Pe^t(\C)$ a point not contained in the
support of $Z$, and $S$ a natural number at most $\deg Z/3$.

\begin{enumerate}

\item
There are fixed constants $c_1,c_2,c_3$ only 
depending on $p$, and $t$, and a subspace $\Pe(F) \subset \Pe^N$ of 
codimension $t-p$, containing $\theta$, intersecting $Z$ properly, and 
fullfilling
\[ c_1 \deg Z \geq D(\Pe(F),Z) \geq - c_2 \deg Z \log \deg Z, \]
and
\[ D^S(Z,\theta) = (D^{\Pe(F)})^S(Z. \Pe(F), \theta) + D(\Pe(F),Z) +
   O((S+\deg Z)\log(S \deg Z)). \]

\item
If $\Pe(F)$ is any subspace of codimension $q \geq t-p$ containing 
$\theta$, then
\[ (D^{\Pe(F)})^S(Z.\Pe(F),\theta) \leq D^S(\theta,Z) + D(\Pe(F),Z) +
    O((S +\deg Z) \log (S \deg Z)). \]
\end{enumerate}
\end{Satz}

The prove will be given for the $D_G$, and $D_{Ch}$ seperately. Consider first
$D_{Ch}$, and
recall that the Chow divisor of an algebraic cycle $X \in Z^p_{eff}(\Pe^t)$ 
is defined in the following way.
Let $\delta: \check{\Pe}^t \to (\check{\Pe}^t)^p$ be the diagonal, and
define the correspondence

\special{em:linewidth 0.4pt}
\linethickness{0.4pt} 
\vspace{4mm}
\hspace{5cm}
\begin{picture}(160.00,40.00)
\put(68.00,28.00){\vector(-3,-2){24.00}}
\put(80.00,35.00){\makebox(0,0)[cc]{$\CC$}}
\put(30.00,5.00){\makebox(0,0)[cc]{$(\Pe^t)^p$}}
\put(130.00,5.00){\makebox(0,0)[cc]{$(\check{\Pe}^t)^p$}}
\put(54,28){\makebox(0,0)[cc]{$f$}}
\put(106,28){\makebox(0,0)[cc]{$g$}}
\put(92.00,28.00){\vector(3,-2){24.00}}
\end{picture}
\vspace{-11mm}
\begin{equation} \label{chow} \end{equation}
\vspace{-2mm}

where $\CC$ is the subscheme of $(\Pe^t)^p \times (\check{\Pe^t})^p$ 
assigning to each $t+1$ dimensional vector space $V$ over a field $k$ the set  
\[ \{ (v_1, \ldots, v_p,\check{v}_1,\ldots \check{v}_p | v_i \in V,\;
   \check{v}_i \in \check{V}, \; \check{v}_i(v_i) = 0, \; 
   \forall i = 1,\ldots p \}. \] 
The maps
$f:\CC \to (\Pe^t)^p, g:\CC \to (\check{\Pe}^t)^p$ are just the restrictions
of the projections. They are flat, projective, surjective,
and smooth.
For $\CZ \in Z^{t+1-p}_{eff}(\Pe^t)$,
the Chow divisor $Ch(\CZ) \subset (\check{\Pe}^t)^p$ is defined as
$Ch(\CZ) := g_* \circ f^* \circ \delta_* (\CZ)$.

\satz{Lemma} \label{PF}
Let $Z \in Z^p_{eff}(\Pe^t_\C)$. 
\begin{enumerate}

\item 
For every $l \geq t+1-p$, there is a subspace
$\Pe(V) \subset \Pe^t(\C)$ of codimension $l$
such that with $V$ the corresponding point in $G_{t+1-l,t+1}$
\[ \log |V, V_Z| \geq -c_1 - \log \deg Z, \]
where $c_1$ is a positive constant only depening on $t$ and $p$.

\item
Let $\Pe(F) \subset \Pe^t(\C)$ be a subspace of codimension $r \leq t+1-p$
that contains a subspace $\Pe(V)$ of codimension $l \geq t+1-p$ with
\[ \log |V, V_Z| \geq - c_1 -\log \deg Z. \]
Then,
\[ D(\Pe(F),Z) \geq - c_3 \deg Z \log \deg Z \]
with $c_3$ a positive constant depending only on $p,q$ and $t$.
\end{enumerate}
\end{Satz}

\proof 
1.
Let $c_2$ be the contstant form Lemma \ref{tub}.2, and
$U(V_Z)$ the tubular neighbourhood of diamater
$\frac1{2 c_2 \deg z}$ of the support of $V_Z$ in $G$. By Lemma \ref{tub}.2,
\[ \mu(U(V_Z)) \leq \frac{c_2 \deg Z}{2c_2 \deg Z} = 1/2 < 1
   = \mu(G). \]
Hence, there is a subspace $\Pe(V)$ of dimension $t-q$ such that the
point $V \in G$ corresponding to $\Pe(V)$ does not lie in
$U(V_Z)$, i.\@ e.\@ $\log |V,V_Z| \geq - \log c_2 -\log (2 \deg Z)$. 
Take $c_1 = 2 c_2$.


\vspace{2mm}

2.
Follows in the same way as \cite{App1}, Proposition 4.17.



\satz{Proposition} \label{raumabl}
For $p+r = t$, let $Z \in Z^p_{eff}(\Pe^t_\C)$, and 
$\Pe(F) \subset \Pe^t(\C)$ be a subspace of codimension $r$,
regular with respect to
$Z$ that contains a subspace $\Pe(V)$ of codimension $l \geq t+1-p$ such that 
in the Grassmanian $G_{p,t+1}$ one has
$\log |V,V_Z| \geq - \log c_1 -\log \deg Z$ with $c_1 > 0$ a constant only 
depending on $p,q$, and $t$. Then, for every $S \leq \deg Z$
\[ D^S(\Pe(F),Z) \leq c_3 ((\deg Z + S)(\log \deg Z + \log S)), \]
and with $U_V$ a neighbourhood of $V$ in $G$, and
$\varphi: U \to U_V$ an affine chart, the function
\[ D_*^S (\Pe(F),X) := \log \sup_{s \leq S, |I|=s}
   \partial^I (\exp - D(\Pe(F),Z)), \]
also fullfills
\[ D_*^S(\Pe(F),X) \leq c_3 ((\deg Z + S) (\log \deg Z+\log S)). \]
with $c_3$ a constant depending only on $p,q,r$, and $t$.
\end{Satz}

\proof
Let $U_F$ be the ball with logarithmic radius $-\log  2c_1 -\log \deg Z$
around $F$ in $G_{t+1-r,t+1}$. By Lemma \ref{trfunk0}, for
every $\tilde{F} \in U_F$ there is a $\tilde{V} \subset \tilde{F}$
of dimension $t-r$
such that $|V,\tilde{V}| \leq |F, \tilde{F}| \leq -2 c_1 - \log \deg Z$.
The assumption on $\Pe(V)$, together with the triangle inequality,
implies $|\tilde{V},V_Z| \geq - \log 2 c_1 - \log \deg Z$, which in turn by
Lemma \ref{PF}.2 implies
\[ D(\Pe(\tilde{F}),Z) \geq -  c_2 \deg Z \log \deg Z, \]
or
\[ D_*(\Pe(\tilde{F}),Z) \leq 2 c_2 \deg Z \log \deg Z. \]
Further by Proposition \ref{prohol}.1 the function
$D(\Pe(\tilde{F}),Z)$ and thereby the function
$D_* (\Pe(\tilde{F},Z))$ has a holomorphic model $g$ on $U_F$. Hence
Proposition \ref{Cauchy}, together with the above inequalities implies
\[ D^S_*(\Pe(F),Z) \leq \sup_{\tilde{F} \in U_F} \log |g(\tilde{F})|
   + 2 c_1 S \log \deg Z + O((S+ \deg Z) \log (S\deg Z))  \leq \]
\[ 2 c_2 \deg Z \log \deg Z + 2 c_1 S \log \deg Z + O(S \log S) \leq
   c_3 (S+ \deg Z) (\log S + \log \deg Z), \]
with a suitable constant $c_3$.

The inequality 
\[ D^S(\Pe(F),X) \leq c_3 ((\deg X + S)(\log \deg X + \log S)) \]
follows in the same way, this time using
$D(\Pe( \tilde{F}),Z) \leq c_4 \deg Z$ for every $\tilde{F}$ regular with 
respect to $Z$ which is just a reformulation of 
\cite{BGS}, Propostions 5.1, and the holomorphic model for
$D(\tilde{F},Z)$.



\satz{Proposition} \label{punktraum}
Let $p,q,r$ be numbers fullfilling $q \geq t+1-p$, $r =t-p$, 
$t-p-q+1 \leq 0$  and
$Z$ be an effecitve cycle of codimension $p$ in $\Pe^t_\C$. Further,
$\Pe(W) \subset \Pe^t(\C)$ a subspace of codimension $q$ that
does not meet the support of $Z$. There is a subspace 
$\Pe(F) \subset \Pe^t$ of codimension $r$ that contains $\Pe(W)$, hence
intersects $Z$ properly such that with $Z_F$ as defined in section 2,
\[ D_\bullet^S(Z,\Pe(W)) \leq D_\bullet^S (Z_F,\Pe(W)) +
   O((\deg Z +S) \log (S\deg Z)), \]
and
\[ D_\bullet^S(Z_F, \Pe(W)) \leq D_\bullet^S(Z,\Pe(F)) +
   O((\deg Z +S) \log (S\deg Z)), \]
where $D_\bullet$ may be chosen to mean either $D_{Ch}$ or $D_G$.
Further, if $\Pe(F)$ is a subspace of codimension $l \geq q$, contains
$\theta$ as well as a subspace $\Pe(V)$ of codimension $t+1-p$ such that 
$\log |\Pe(V),Z| \geq -c -\log \deg Z$, and $Z . \Pe(\tilde{F})$ for
every $\tilde{F}$ in some neighbbourhood of $F$ is a sum
of projective subspaces, then the above inequalities still hold.

\end{Satz}

\proof 
Let $\Pe(V) \subset \Pe(F)$ be a subspace of codimension 
$l = 2t+1-p-q \geq t+1-p$ such that $V$ has
maximal distance to the support of $V_Z$ with this property.
By Lemma \ref{PF}, $\log |V,V_Z| \geq -c_1 \log \deg Z$. Let further $\Pe(F)$
be the unique subspace of $\Pe^t$ that contains $\Pe(W)$ as well as $\Pe(V)$,
and $U_W$ a neighbourhood of $W$ in $G_{t+1-q,t}$ such that for
each $\tilde{W} \in U_W$ the intersection of $\Pe(\tilde{W})$ with $Z$ is 
proper. Finally, let $f: U_W \to G_{t+1-q,t}$ be the map from Lemma
\ref{trfunk}. By Proposition \ref{raumraum} for every $\tilde{W} \in U_W$
and $\tilde{F} = f(\tilde{W})$,
\begin{equation} \label{punktraum1}
D(Z_{\tilde{F}},\Pe(\tilde{W})) - D(Z,\Pe(\tilde{W})) = 
D(Z_{\tilde{F}},\Pe(\tilde{F}))- D(Z,\Pe(\tilde{F})), 
\end{equation}
where $D(Z_{\tilde{F}},\Pe(\tilde{F})) = c \deg Z$ by \cite{App1}, . 
Let
\[ \varphi: \A^{(t+1-r)r}   \to G_{t+1-r,t}, 
   \quad \A^{(t+1-r)r} \to G_{t+1-q,t} \]
be the canonical affine charts from chapter 3 centered at $W$ and $F$ 
respectively, and define
\[ G(\tilde{W},\tilde{F}) = \exp (D(Z_{\tilde{F}}, \Pe(\tilde{W}))) =
   \exp(D(Z_{f(\tilde{W})},\Pe(\tilde{W})), \]
\[ F(\tilde{W},) = \exp (D(Z,\Pe(\tilde{W}))), \quad 
   \mbox{and} \quad H(\tilde{F}) = \exp(D(Z,\Pe(\tilde{F}))). \]
Then, (\ref{punktraum1}) reads
\[ F(\tilde{W}) = \exp(-c \deg Z) \;
   H \circ f (\tilde{W}) \; G(\tilde{W}, f \tilde{W}), \]
hence for every $I$ with $|I|= s \leq S$,
\begin{eqnarray*}
  \partial^I ((\varphi^* F) &=& \exp(- c \deg Z) 
   \partial^I ((H \circ f \circ \varphi) 
   (G \circ (\varphi, f \circ\varphi))) \\
&=& \exp(-c \deg Z)
   \partial^I ((H \circ \psi \circ \psi^{-1} \circ f \circ \varphi) 
   (G \circ \psi \circ \psi^{-1} \circ (\varphi, f \circ \varphi))). 
\end{eqnarray*}
By Lemma \ref{PF}, 
$(\partial^J (\psi^{-1} \circ f \circ \varphi))(0) \leq cS \log S$
for every $J \prec I$, and by the previous Lemma 
\[ |(\partial^I (H \circ \psi))(0)| \leq c_1 (\deg Z + S) \log(S \deg Z), \]
which for every $I$ with $|I| \leq S$, which by elementary differentiation 
techniques implies
\[ |(\partial^I(\varphi^* F))(0)| \leq
   \sup_{|J| \leq S} |(\partial_W^J (\psi^*G))(0)| + c S \log S +
   c_1 (S + \deg Z) \log (S \deg Z) + c_2 S, \]
where $\partial_W^I$ denotes partial derivatives by the first component; 
hence
\[ D^S(Z, \Pe(W)) \leq D^S(Z_F, \Pe(W)) - c \deg Z + c S \log S +
   c_1 (S + \deg Z) \log (S \deg Z) + c_2 S. \]
The inequaltiy in the other direction is proved analogously.

\satz{Lemma} \label{DXF}
Let $Z \subset \Pe^t_\C$ be a subvariety of codimension $p$, and
$\Pe(W) \subset \Pe^t_\C$ a subspace of codimension $q \geq t-p$ that
does not meet $Z$. Finally $\Pe(F)$ a subspace of
codimension $r = t-p$ containing $\Pe(W)$, and intersecting $Z$ properly,
\[ D^S(\Pe(W),Z_F) = (D^{\Pe(F)})^S (\Pe(W), Z_F . \Pe(F)) + 
   c_{\cdots} \deg Z. \]
\end{Satz}

\proof
Let $\varphi: \A^{(t+1-q)q} \to G_{t+1-q,t+1}$ be an affine chart centered
at the origin such that if $G^F \subset G_{t+1-q,t+1}$ is the Grassmannian
of $t+1-q$-dimensional subspaces that are contained in $F$, and
$\A^{(t+1-q)(p+q-t)} \subset \A^{(t+1-q)q}$ is the affine subset
corresponding to the first $(t+1-q)(p+1-t)$  coordinates. The restriction
of $\varphi$ to $\A^{(t+1-q)(p+q-t)}$ is an affine chart for $G^F$.
Further, let $(x_1,y_1, \ldots, x_{(t+1-q)(p+q-t)}, y_{(t+1-q)(p+q-t)})$ be
the real coordinates of $\A^{(t+1-q)(p+q-t)}$, and denote by
$\partial_{x_1}, \ldots, \partial_{y_{(t+1-q)(p+q-t)}}$ or simply
$\partial_1, \ldots, \partial_{2(t+1-q)(p+q-t)}$ the partial derivatives
with respect to these coordinates, and for 
$I_F = (i_1, \ldots, i_{2(t+1-q)(p+q-t)}$ let $\partial_F^{I_F}$ be
the differential operator
$\partial_1^{i_1} \ldots \partial_{2(t+1-q)(p+q-t)}^{i_{2(t+1-q)(p+q-t)}}$.

By Lemma \ref{XF},
\[ D(\Pe(W),Z_F) = D^{\Pe(F)}(\Pe(W),Z_F . \Pe(F)) + c \deg Z. \]
Since, with these notations, for $i = 2(t+1-q)(p+q-t)+1, \ldots, 2(t+1-q)(t+1)$
we have $\partial_i D(Z_F, \Pe(W))= 0$, it follows
\[ D^{\Pe(F)})^S (\Pe(W), Z_F . \Pe(F)) =
   \log \left| \sup_{s \leq S, |I_F|= s} \partial_F^{I_F} 
      \exp D(\Pe(W), Z_F. \Pe(F)) \right| = \]
\[ \left| \sup_{s \leq S, |I|=s} \partial^I 
              \exp (D(\Pe(W),Z_F) - c \deg Z) \right| =
   D^S (\Pe(W),Z_F) - c \deg Z. \]

\satz{Corollary} \label{punktraum2}
In the situation of Proposition \ref{raumabl},
\[ D_\bullet^S(Z,\Pe(W)) \leq (D_\bullet^{\Pe(F)})^S (Z . \Pe(F),\Pe(W)) +
   O((\deg Z +S) \log (S\deg Z)), \]
and
\[ (D_\bullet^{\Pe(F)})^S(Z . \Pe(F), \Pe(W)) \leq D_\bullet^S (Z,\theta) +
   O((\deg Z +S) \log (S\deg Z)), \]
\end{Satz}

\proof
Follows immediately from proposition \ref{raumabl} and the previous Lemma.

\proof {\sc of Theorem \ref{zerl}}
In the Corollary, take $q=t$.

\proof {\sc of Theorem \ref{main2}:}
Let $\Pe(F)$ be as in Theorem \ref{zerl}. Then,
\[ D^S(\theta,Z) \leq (D^{\Pe(F)})^S(\theta, Z. \Pe(F)) + 
   O((S\deg Z) \log (S + \deg Z)). \]
Since 
\[ Z. \Pe(F) = \sum_{i=1}^{\deg Z} z_i \]
with $|z_1,\theta| \leq \cdots \leq |z_{\deg Z},\theta|$ is zero-dimensional, 
by Proposition \ref{Punkte},
\[ (D^{\Pe(F)})^S(\theta,Z. \Pe(F)) \leq 
   \sum_{i=S+1}^{\deg Z} \log |z_i,\theta| + O(S \log \deg Z). \]
The two inequalities together imply
\[ D^S(\theta,Z) \leq \sum_{i=S+1}^{\deg Z} \log |z_i,\theta| +
   O((S + \deg Z) \log (S \deg Z)). \]
Similarly,
\[ (D^{\Pe(F)})^{3S}(\theta, Z. \Pe(F)) \leq D^{3S} (\theta, Z) +
   O((S+\deg Z) \log (S \deg Z), \]
and
\[ 2 \sum_{S+1}^{\deg Z} \log |z_i,\theta| \leq
   (D^{3S})^{\Pe(F)}(\theta,Z . \Pe(F)) + O((S+\deg Z) \log (S\deg Z)), \]
imply
\[ \sum_{S+1}^{\deg Z} 2 \log |z_i,\theta| 
   \leq D^{3S}(\theta,Z) + O(\deg Z \log \deg Z). \]

\section{Proof of the main Theorems}

\subsection{Combinatorics of the intersection points}

Let $Z_0,Z_1$ be properly intersecting effective cycles of
pure codimensions $p$ and $q$ in $\Pe^t(\C)$, and $Z_0 \# Z_1$ their join. 
(See \cite{App1}, section 6 for details.) Let further $\theta \in \Pe^t(\C)$
be a point not contained in the support of $Z_0 . Z_1$,
and $\Pe(F_0), \Pe(F_1) \subset \Pe^t$ projective subspaces of
dimensions $p,q$ such that $\Pe(F_i)$ intersects $Z_i$ properly
for $i=0,1$. Denote
\[ Z_i . \Pe(F_i) = \sum_{j=1}^{\deg Z_i} z^i_j \quad \mbox{with} \quad
   |z^i_1,\theta| \leq \cdots \leq |z^i_{\deg Z_j},\theta|, \quad j=0,1, \]
and assume that the numbers
\[ |z^0_1,\theta|, \ldots, |z^0_{\deg Z_0},\theta|, 
   |z^1_1,\theta|, \ldots, |z^1_{\deg Z_1},\theta| \]
together with the numbers
\[ |z^0_i \# z^1_j,(\theta,\theta)|, \quad i=1, \ldots  \deg Z_0,
   \quad j=1, \ldots,Z_1 \]
are pairwise distinct.

Let next $\underline{n}$ for $n \in \N$ denote the set
$\{1,\ldots,n\}$ and define a path
\[ f = f_\theta = (f_0,f_1): \underline{\deg Z_0 + Z_1} 
   \to \underline{\deg Z_0} \times \underline{\deg Z_1} \]
as in the introduction 
in the following way: $f(k) = (f_0(k),f_1(k)) = (\nu_0,\nu_1)$, 
iff $k=\nu_0 + \nu_1$ and there is a number $t \in [0,1]$ such that
\[ |z^0_\nu,\theta| < t < |z^0_{\nu+1},\theta| \quad \mbox{and} \quad
   |z^1_\kappa,\theta| < t < |z^1_{\kappa+1},\theta|. \]
The maps $f_0$ and $f_1$ are surjective. Further for 
$k \in \underline{\deg Z_0 + \deg Z_1}$, set $i_k = 0$ 
if $f_1(k) > f_1(k-1)$, hence 
$|z^1_{f_1(k)},\theta|<|z^0_{f_0(k-1)-1},\theta|$, and $i_k = 1$ otherwise.

Recall that for $x,y,\theta \in \Pe^t(\C)$ the inequalities
\begin{equation} \label{jpabsch}
\mbox{min}(|x,\theta|,|y,\theta|) \leq |x \# y,(\theta,\theta)|
\leq \mbox{max}(|x,\theta|,|y,\theta|) 
\end{equation}
hold. (\cite{App1}, Lemma 6.4)

\satz{Lemma} \label{comb1}

\begin{enumerate}

\item
With the above notations, let $K$ be a number such that either
$f_0(K-1) < \deg Z_0$ or $f_1(K_1) < \deg Z_1$. Then,

\[ \sum_{i=1}^{\deg Z_0} \sum_{j=1}^{\deg Z_1}
   \log |(\theta,\theta),z^0_i \# z^1_j| \leq \;
   \sum_{k=0}^K \sum_{l=f_{i_k}(k)+1}^{\deg Z_{i_k}} 
   \log |\theta,z^{i_k}_l|. \]

\item
For every $K \leq \deg Z_0 + \deg Z_1$, and $(\nu_1,\nu_0)= f(K)$,
\[ \sum_{i=1}^{\deg X_0} \sum_{j=1}^{\deg X_1}
   \log |(\theta,\theta),z^0_i \# z^1_j| \leq   
   \nu_1 \sum_{k= \nu_0+1}^{\deg Z_0} \log |\theta, z_k^0| +
   \nu_0 \sum_{k= \nu_1}+1^{\deg Z_1} \log |\theta, z_k^1|. \]
\end{enumerate}
\end{Satz}

\proof
1.
The equation
\[ \sum_{i=1}^{\deg Z_0} \sum_{j=1}^{\deg Z_1} 
   \log |(\theta,\theta),z_i^0 \# z_j^1| =
   \sum_{k=1}^K \; \sum_{l= f_{i_k}(k)+1}^{\deg Z_{i_k}}
   \log |(\theta,\theta),z^{1-i_k}_{f_{i_k}(k)} \# z^{i_k}_l| \]
is just a reordering of the sum.
By (\ref{jpabsch}), and the fact 
$|\theta, z^{1-i_k}_{f_{i_k}(k)}| \leq |\theta,z^{i_k}_l|$ following 
immediately from the definition of $f$, and $i_k$, the left hand side is less 
or equal
\[ \sum_{k=1}^K \sum_{l=f_{i_k}(k)+1}^{\deg Z_{i_k}} 
   \log |\theta,z^{i_k}_l|, \]
as was to be proved.

\vspace{2mm}

2.
Follows from 1, by taking on the right hand side only the first
$\nu_0 + \nu_1 \leq K$ summands in the first summation, and the
last $\deg Z_{i_k} - \nu_{i_k}$ summands in the second summation, which
is possible since $\nu_{i_k} \geq f_{i_k}(k)$.

\vspace{2mm}

Let now $S \leq \deg Z_0 \deg Z_1$, and 
$M \subset \{z^0_i \# z^1_j | i=1, \ldots, \deg Z_0, j=1, \ldots, \deg Z_1\}$
be the subset with $|M|=S$, and $|z,(\theta,\theta)| < |z',(\theta,\theta)|$ 
for every $z \in M, z' \notin M$. 
Because by (\ref{jpabsch}), 
\[ |x,\theta| \leq |x',\theta| \quad \mbox{and} \quad
   |y,\theta| \leq |y',\theta|, \quad \Longrightarrow \quad
   |x \# y,(\theta,\theta)| \leq |x' \#  y',(\theta,\theta)| \]
the set $M$ fullfills the condition 
\begin{equation} \label{ablpbed} 
\forall i,i' \in \underline{\deg Z_0}, \quad 
   \forall j,j' \in \underline{Z_1} \quad i \leq i' \wedge j \leq j' \wedge
   (i,j) \notin M \Longrightarrow (i',j') \notin M;
\end{equation}
consequently, there is a number $k_0 \leq K$ such that
$f(k_0) = (\nu_0,\nu_1) \notin M$, and $\nu_0 \nu_1 \leq S$.
For any $k \geq k_0$ also $f_0(k) \geq \nu_0$, $f_1(k) \geq \nu_1$,
and $f(k) \notin M$.
Hence, if $h_S(k)$ denotes the number 
$\mbox{min}\{ l | z^{1-i_k}_{f_{i_k}(k)} \# z_l^{i_k} \notin M \}-f_{i_k}$,
we have $h_S(k) = 0$ for $k \geq k_0$.

\satz{Lemma} \label{comb2}
With the above notations,
\begin{enumerate}
\item
\[ \sum_{z \notin M} \log |(\theta,\theta),z| \leq
   \sum_{k=k_0}^K \; \sum_{l=f_{i_k}(k)-h_S(k)+1}^{\deg Z_{i_k}}
   \log |\theta,z^{i_k}_l|. \]

\item
For any $k \geq k_0$, and  $(\bar{\nu}_0,\bar{\nu}_1) = f(k)$,
\[ \sum_{z \notin M} \log |(\theta,\theta),z| \leq
   (\bar{\nu}_1 -\nu_1) \sum_{l= \bar{\nu}_0+1}^{\deg Z_0} 
   \log |\theta, z_l^0| +
   (\bar{\nu}_0-\nu_0) \sum_{l= \bar{\nu}_1´1}^{\deg Z_1} 
   \log |\theta, z_l^1|. \]

\item
With $k, \bar{\nu}_0,\bar{\nu}_1$ as in 2,
\[ (\bar{\nu}_0-\nu_0) (\bar{\nu}_1-\nu_1) 
   \log |Z_0+Z_1,\theta| +
   \sum_{z \notin M} \log |(\theta,\theta),z| \leq \]
\[ (\bar{\nu}_1-\nu_1) \sum_{l=\nu_0+1}^{\deg Z_0}  \log |\theta,z^0_l| +
   \bar{\nu}_0-\nu_0) \sum_{l=\nu_1+1}^{\deg Z_1} \log |\theta,z^1_l|. \]

\end{enumerate}
\end{Satz}

\proof
1. 
Since $z^0_i \# z^1_j \notin M$ implies $i \geq \nu_0$ or $j \geq \nu_1$, the
inequality
\[ \sum_{z \notin M} \log |(\theta,\theta),z| \leq
   \sum_{k = k_0}^K \; \sum_{l=f_{i_k}(k)-h_S(k)+1}^{\deg Z_{i_k}}
   \log |(\theta,\theta),z^{1-i_k}_{f_{i_k}(k)} \# z_l^{i_k}| \]
again follows from a renumbering of the sum.
The inequality
\[ \sum_{k = k_0}^K \sum_{l=f_{i_k}(k)}^{\deg Z_{i_k}}
   \log |(\theta,\theta),z^{1-i_k}_{f_{i_k}(k)} \# z_l^{i_k}| \leq
   \sum_{k=k_0}^K \sum_{f_{i_k}+1}^{\deg Z_{i_k}}
   \log |\theta,z^{i_k}_l| \]
follows from (\ref{jpabsch}) as in the previous Lemma.

\vspace{2mm}

2.
Follows again from 1 by leaving out the first 
$\bar{\nu}_1+ \bar{\nu}_0- \nu_1-\nu_0$ in the first summation,
and taking only the last 
$\deg Z_{i_k} - \bar{\nu}$ summands in the second summation.

\vspace{2mm}

3. 
Define the sets
\[ N_0 := \{ (i,j)  \in \underline{\deg Z_0} \times \underline{\deg Z_1} |
   \nu_0 \leq i \leq \bar{\nu}_0 \}, \]
\[ N_1 := \{ (i,j)  \in \underline{\deg Z_0} \times \underline{\deg Z_1} |
   \nu_1 \leq j \leq \bar{\nu}_1 \}. \]
The set $N_0 \cap N_1$ has cardinality 
$(\bar{\nu}_0-\nu_0)(\bar{\nu}_1-\nu_1)$.
Thus the first inequaltity of (\ref{jpabsch}) implies
\[ (\bar{\nu}_0-\nu_0)(\bar{\nu}_1-\nu_1) \log |Z_0 + Z_1, \theta| \leq
   \sum_{z \in N_0 \cap N_1} \log |(\theta,\theta),z| \]
for each $z \in N_0 \cap N_1$

Further, $N_0 \cap N_1$ is contained in the 
complement of $M$. Hence, by the second inequality of (\ref{jpabsch}),
\[ (\bar{\nu}_0-\nu_0) (\bar{\nu}_1-\nu_1) 
   \log |Z_0 + Z_1, \theta| + \sum_{z \in N_0 \cap N_1}
   \log |(\theta,\theta),z| \leq \]
\[ \sum_{l=\nu_0+1}^{\bar{\nu}_0} |(\theta,\theta),z^0_l| +
    \sum_{l=\nu_1+1}^{\bar{\nu}_1} |(\theta,\theta),z^1_l|. \]
Adding the equality from part 2 of the Lemma proves the claim.

\subsection{Finish of proofs}

By Lemma \ref{PF} and Proposition \ref{raumabl}, 
the $\Pe(F_i)$ $i=0,1$ from the previous section may be chosen in such a way 
that they contain subspaces $\Pe(V_i)$ of codimension
$t+1- \mbox{codim} \; Z_i$ such that 
\[ \log |\Pe(V_i),Z_i| \geq -c -\deg Z_i, \]
hence
\[ D^S(\Pe(F_i),Z_i) \geq - O((\deg Z_i +S) (\log \deg Z_i + log S). \]
for every $S$, and by Theorem \ref{zerl}, 
\[ 2\sum_{j=S+1}^{\deg Z_0} \log |z^0_j,\theta| \leq D^{3S}(\theta,Z_0) +
   O((S+\deg Z_0) \log (S \deg Z_0)), \]
\begin{equation} \label{main1eins}
2\sum_{j=\bar{S}+1}^{\deg Z_1} \log |z_j^1,\theta| \leq D^{3S}(\theta,Z_1)
+ O((S+\deg Z_1) \log \deg (S Z_1)) 
\end{equation}
Further, we have
\[ (\Pe(F_0) \# \Pe(F_1)) . (Z_0 \# Z_1) =
   \sum_{i=1}^{\deg Z_0} \sum_{j=1}^{\deg Z_1} z_i^0 \# z_j^1, \]
and each of the $z_i^0 \# z_j^1$ is a one dimensional
projective subspace of $\Pe^{2t+1}$.
We denote
\begin{equation} \label{ach}
(\Pe(F_0) \# \Pe(F_1)) . (Z_0 \# Z_1) = \sum_{i=1}^{\deg Z_0 \deg Z_1} z_i, 
\end{equation}
such that $|(\theta,\theta),z_1| < \ldots < 
|(\theta,\theta),z_{\deg Z_0 \deg Z_1}|$.


\satz{Proposition} \label{eins}
With the above notations, and $K$ as in Lemma \ref{comb1}
\begin{enumerate}

\item
\[ 2 \sum_{i=1}^{\deg Z_0} \sum_{j=1}^{\deg Z_1} 
   \log |(\theta,\theta), z^0_i \# z^1_j| 
   \leq \sum_{k=0}^K D^{3(f_{i_k}(k)-h_S(k))} (Z_{i_k},\theta) + \]
\[ O(\deg Z_0 \deg Z_1 \log(\deg Z_0 \deg Z_1)). \]

\item
With $S,M, \nu_0, \nu_1, k_0$ as in Lemma \ref{comb2}, 
\[ 2\sum_{z \notin M} \log |(\theta,\theta), z| 
   \leq \sum_{k=k_0}^K D^{f_{i_k}} (Z_{i_K},\theta) + \]
\[ O(\deg Z_0 \deg Z_1 \log(\deg Z_0 \deg Z_1)). \]

\item
For every $k \geq k_0$, and $(\bar{\nu}_0,\bar{\nu}_1) = f(k)$,
\[ 2(\nu-\nu_0) (\kappa-\kappa_0) \log |Z_0+Z_2,\theta|+
   2\sum_{z \notin M} \log |(\theta,\theta), z| \leq \]
\[ (\bar{\nu}_1-\nu_1) D^{\nu_0}(\theta,Z_0) + 
   (\bar{\nu}_1- \nu_1) D^{\nu_1}(\theta,Z_1)+
   O(\deg Z_0 \deg Z_1 \log(\deg Z_0 \deg Z_1)). \]
\end{enumerate}
\end{Satz}

\proof
1.
By Lemma \ref{comb1}.1, 
\begin{equation} \label{main1zwei}
2 \sum_{i=1}^{\deg Z_0} \sum_{j=1}^{\deg Z_1} 
   \log |(\theta,\theta), z^0_i \# z^1_j| 
   \leq \sum_{k=0}^K \sum_{l=f_{i_k}(k)-h_S(k)+1}^{\deg Z_{i_k}} 
   \log |\theta,z^{i_k}_l|, 
\end{equation}
Now for each $k$,
\[ 2 \sum_{l=f_{i_k}(k)-h_S(k)+1}^{\deg Z_{i_k}} \log |\theta,z^{i_k}_l| \leq
   D^{3(f_{i_k}(k)-h_S(k)} (Z_{i_k},\theta) + 
   c\deg Z_{i_k} \log \deg Z_{i_k}, \]
by (\ref{main1eins}). Since $i_k = 0$ for at most $\deg Z_1$ values
of $k$ and $i_k=1$ for at most $\deg Z_0$ values of $k$, we have
\[ \sum_{k=0}^K \deg Z_{i_k} \log \deg Z_{i_k} \leq
   2c\deg Z_0 \deg Z_1 \log (\deg Z_0 \deg Z_1). \]
Hence, the left hand side of (\ref{main1zwei}) is less or equal
\[ \sum_{k=1}^K D^{f_{i_k}(k)} (Z_{i_k},\theta) +
   2\deg Z_0 \deg Z_1 \log (\deg Z_0 \deg Z_1), \]
and the claim follows.

\vspace{2mm}

2. 
Follows in exactly the same way as 1, this time using Lemma \ref{comb2}.1.

\vspace{2mm}

3.Follows from Lemma \ref{comb2}.3.




\satz{Lemma} 
Let $\Pe(V) \subset \Pe(F)$ be a subspace of codimension $t-p+1$
such that $\log |\Pe(V),X| \geq - c \log \deg X$, and
$\Pe(V') \subset \Pe(F')$ a subspace of codimension $t-q+1$ such that
$\log |\Pe(V'),Y| \geq -c \log \deg Y$. Then,
\[ \log |\Pe(V) \# \Pe(V'), X \# Y| \geq - c \log (\deg X \deg Y). \]
\end{Satz}

\proof
Follows from the inequality
\[ |x \# y, v \# w| \geq \mbox{min}(|x,v|, |y,w|) \]
from \cite{App1}, Lemma 6.4.

\vspace{2mm}

By this Lemma the pair $\Pe(F_0) \# \Pe(F_1),Z_0 \# Z_1$ fullfills
the condition of Proposition \ref{raumabl}.

\satz{Proposition} \label{zwei}
With the above notaions from (\ref{ach}),
\[ D^S(\theta,Z_0.Z_1) + D(Z_0,Z_1) \leq \]
\[ \sum_{i=S+1}^{\deg Z_0 \deg Z_1} \log |(\theta,\theta), z_i|+
   O((S+\deg Z_0 \deg Z_1) \log (s \deg Z_0 \deg Z_1). \]
\end{Satz}

\proof
Firstly, by Proposition \ref{Punkte}, 
\[ D^S((\theta,\theta), X\#Y) \leq
   \sum_{i=S+1}^{\deg Z_0 \deg Z_1} \log |(\theta,\theta),z_i| + \]
\[ O(S \log(\deg Z_0 \deg Z_1)). \]
Next, together with the preceeding Lemma, the
Propositions \ref{punktraum}, and \ref{DXF}  just as in the proof of
Corollary \ref{punktraum2} imply
\begin{eqnarray*}
(D^{\Pe(\Delta)})^S ((\theta,\theta), (Z_0 \# Z_1) . \Pe(\Delta)) &\leq&
   D^S((\theta,\theta), Z_0\#Z_1) - D(\Pe(\Delta),(X \# Y)) + \\ 
&& O((S+\deg Z_0 \deg Z_1) \log (S \deg Z_0 \deg Z_1)). 
\end{eqnarray*}
Since $(X \# Y) . \Pe(\Delta)) = \delta(X. Y)$
with $\delta: \Pe^t \to \Pe(\Delta) \subset \Pe^{2t+1}$ the diagonal map,
\[ D^{\Pe(\Delta)} ((\theta,\theta),(Z_0 \# Z_1) . \Pe(\Delta)) =
   D(\theta, X . Y) + \log 2 \deg Z_0 \deg Z_1, \]
and consequently
\[ (D^{\Pe(\Delta)})^S ((\theta,\theta),(X \# Y) . \Pe(\Delta)) =
   D^S(\theta, X . Y) + \log 2 \deg Z_0 \deg Z_1. \]
Similarly,
\[ D(\Pe(\Delta), X\#Y) = D(X,Y) + \log 2 \deg Z_0 \deg Z_1. \]
The claim follows.

\proof {\sc of Theorem \ref{DMBT1}}
Follows from Proposition \ref{eins}.1 together with Proposition
\ref{zwei} for $S=0$.

\proof {\sc of Corollary \ref{cor1}:}
Since $f_{i_k}(k) \leq f_{i_k}(l)$ for $k \leq l$, we get 
\[ D^{f_{i_k}(k)} (Z_{i_k},\theta) \leq
   D^{f_{i_k}(l)} (Z_{i_k},\theta) = D^{\nu_i}(Z_{i_k},\theta), \]
and the claim follows from Theorem \ref{DMBT1} by cutting the sum on the left
hand side at $l$.

\proof {\sc of Corollary \ref{cor2}:}
Since in the path $f$ in each step exacly one coordinate increase, there is a
$k \leq S$ such that either $f(k) = (1,\nu_1)$ with $\nu_1 \leq S$
or $f(k) = (0,S)$. In the first case, by Corollary \ref{cor1},
\[ 2D(Z_0,Z_1) + 2 D(Z_0.Z_1,\theta) \leq D^{3 \nu_1} (Z_1,\theta) +
   O((S + \deg Z_0 \deg Z_1) \log(S \deg Z_0 \deg Z_1), \]
wich trivially is less or equal
\[ D^{3 S} (Z_0,\theta) +
   O((S + \deg Z_0 \deg Z_1) \log(S \deg Z_0 \deg Z_1). \]
In the second case, by the same Corollary,
\[ 2 D(Z_0,Z_1) + 2 D(Z_0.Z_1,\theta) \leq S D (Z_0,\theta) +
   O((S + \deg Z_0 \deg Z_1) \log(S \deg Z_0 \deg Z_1). \]

\proof {\sc of Theorem \ref{DMBT2}}
Follows form the Propositions \ref{eins}, and \ref{zwei}.

\proof {\sc of Corollary \ref{cor3}:}
Follows in the same way as Corollary \ref{cor1}.

\proof {\sc of Corollary \ref{cor4}:}
Let $k_0, \nu_1, \nu_1$ be as abouve. We have $\nu_0 \nu_1 \leq S = S_0 S_1$,
and without loss of generality one may assume $S_0 \geq \nu_0 = f_0(k_0)$.
Let $l \geq k_0$ be the smallest number such that $f_0(l) = 2 S_0$,
and $(\bar{\nu}_0,\bar{\nu}_1) = f(l)$. If $\bar{\nu_1}-\nu_1 \geq S_1$
then, by Corollary \ref{cor3}, with $k=l$, 
\begin{eqnarray*}
   2D(Z_0,Z_1) + 2 D^S(Z_0.Z_1,\theta) &\leq& (\bar{\nu_1}-\nu_1)
   D^{6 S_0} (Z_0,\theta) + \\
   &&O((S + \deg Z_0 \deg Z_1) \log(S \deg Z_0 \deg Z_1) \\ \\
   &\leq&
   S_1 D^{6S_0} (Z_0,\theta) + \\
   && O((S + \deg Z_0 \deg Z_1) \log(S \deg Z_0 \deg Z_1). 
\end{eqnarray*}
If $\bar{\nu}_1-\nu_1 \leq S_1$, and $\nu_1 \leq 2 S_1$, then
$\bar{\nu}_1 \leq 3 S_1$, hence
\begin{eqnarray*}
   2D(Z_0,Z_1) + 2 D^S(Z_0.Z_1,\theta) &\leq& (\bar{\nu}_0-\nu_0)
   D^{3 \bar{\nu}_1} (Z_0,\theta) + \\
   && O((S + \deg Z_0 \deg Z_1) \log(S \deg Z_0 \deg Z_1) \\ \\
   &\leq& S_0 D^{9 S_1} (Z_0,\theta) + \\
   && O((S + \deg Z_0 \deg Z_1) \log(S \deg Z_0 \deg Z_1), 
\end{eqnarray*}
since $\bar{\nu}_0-\nu_0 \geq 2S_0-S_0=S_0$, also by Corollary \ref{cor3}.

Finally if $\bar{\nu}_1-\nu_1 \leq S_1$, and $\nu_1 \geq 2 S_1$, then
$\nu_1/2 S_0 \geq S$, hence the complement
of the set $M$ from Lemma \ref{comb2} is contained in the set
\[ \{z_i^0 \#´z_j^1 | i \geq \nu_1/2 \wedge j \geq S_0 \}. \]
This means that for $\nu_1/2 \leq k \leq \nu_1$ and $i_k=0$
the value $h_S(k)$ is less or equal $f_0(k)-S_0$.
Hence, by Theorem \ref{DMBT2}.1, with $l \leq k_0$ the smallest number
such that $f_1(l) = \nu_1/2$
\ref{eins}, and \ref{zwei},  we get 
\begin{eqnarray*}
   2D(Z_0,Z_1) +D^S(Z_0,Z_1,\theta) &\leq& 
   \sum_{k=l+1}^{k_0} D^{3(f_{i_k}-h_S(k))(Z_{i_k},\theta} + \\
   && O((S + \deg Z_0 \deg Z_1) \log(S \deg Z_0 \deg Z_1) \\ \\
   &\leq&
   \sum_{k=\nu_1/2+1}^{\nu_1} D^{3S_0}(Z_0,\theta) + \\
   && O((S + \deg Z_0 \deg Z_1) \log(S \deg Z_0 \deg Z_1) \\ \\ 
   &\leq& S_1 D^{3S_0} (Z_0,\theta) + \\
   && O((S + \deg Z_0 \deg Z_1) \log(S \deg Z_0 \deg Z_1).
\end{eqnarray*}

\begin{appendix}

\section{Proof of Lemmas \ref{hilf}, \ref{hilfzwei} and \ref{hilfdrei}}

\proof{\sc or Lemma \ref{hilf}}
Firstly,
\[ f^{(s)}(0) = s! \sum_{i_1=1}^n \sum_{i_2=i_1+1}^n \cdots 
   \sum_{i_s=i_{s-1}+1}^n \frac{f(0)}{\prod_{k=1}^s x_{i_k}}, \]
for every $s \leq n$. Since 
$\prod_{k=1}^s |x_{i_k}| \geq \prod_{k=1}^s |x_k|$ for every
$s$-tupel $(x_{i_1},\ldots,x_{i_k})$ conseqently,
\begin{equation} \label{hilf0}
|f^{(s)}(0)| \leq \frac{n!}{(n-s)!} \frac{|f(0)|}{\prod_{k=1}^s |x_k|} =
\frac{n!}{(n-s)!} \prod_{i=s+1}^n |x_i|, 
\end{equation}
proving the second inequality.

\vspace{2mm}

The first inequality
will be proved for $|x_1|<|x_2|< \cdots <|x_n|$, and follows for
$|x_1| \leq \cdots \leq |x_n|$ by continuity.

\vspace{2mm}

{\sc Step 1:}
There are points 
\[ x_{11}, \ldots, x_{1n}, x_{22}, \ldots, x_{2n}, \ldots,
   x_{n-1 n-1}, x_{n-1 n}, x_{nn} \in [-1,1], \]
such that
\[ 0<|x_{ii}|< |x_{i i+1}| < \cdots < |x_{in}| < 1, \quad
   \sgn(x_{ij}) = \sgn(x_{i-1j}), \quad \forall i=1,\ldots,n,  \]
\[ 0<|x_{ij}| < |x_{i-1j}| \leq 1, \quad i=2,\ldots,n, j=i, \ldots,n, \]
\[ f^{(s-1)}(x_{sj}) = 0, \quad s=1, \ldots, n, j=s, \ldots,n. \]

\proof
The points are defined recursively. With
$x_{1j} = x_j$ the claims are fullfilled or empty.
Assume the points $x_{ss}, \ldots, x_{sn}$ are defined.
If $\sgn x_{sj} = \sgn x_{s j+1}$, since 
$f^{(s-1)}(x_{sj}) = f^{(s-1)}(x_{sj+1}) =0$, there is a point
$x_{s+1 j+1} \in (x_{sj},x_{sj+1})$, if $\sgn x_{sj}=1$, and
$x_{s+1 j+1} \in (x_{sj+1},x_{sj})$ if $\sgn x_{sj}=-1$ such that
$f^{(s)} (x_{s+1 j+1}) = 0$. If $\sgn x_{sj} \neq x_{sj+1}$, let
$k$ be the biggest number less than $j$ such that 
$\sgn x_{sk} = \sgn x_{s j+1}$; if there is no such number
take $0$ instead of $x_{sk}$ in the following step. Then there is
a $x_{s+1 j+1} \in (x_{sk},x_{s j+1})$ if $\sgn x_{sk}=1$, and
$x_{s+1 j+1} \in (x_{s j+1},x_{sk})$ if $\sgn x_{sk}=-1$ such that
$f^{(s)}(x_{s+1 j+1} =0$. The so constructed points obviously fullfill
the required conditions.

\vspace{2mm}

{\sc Step 2:} 
For every $s$ with $0 \leq s < n$ there is a point $\bar{x}_s \in [-1,1]$
with $|\bar{x}_s| < |x_s|$, and
\begin{equation} \label{hilf1}
|f^{(s)}(\bar{x}_s)| \geq \frac{1}{2^s} \prod_{i=s+1}^n |x_i| > 0. 
\end{equation}

\proof
The claim obviously holds for $s=0$ with $\bar{x}_0 = 0$. 
Let $y_i = x_{i+1i+1}$ for $i=1, \ldots, n$ with $x_{i+1i+1}$ the points from 
step 1, and assume the claim holds for $s$. Then, by (\ref{hilf1}), and
step 1, $f^{(s)} (\bar{x}_s) \neq 0$, and $f^{(s)}(y_s) = 0$, we have 
$\bar{x}_s \neq y_s$. By the mean value Theorem, there is a $\bar{x}_{s+1}$
insinde the intervall between $\bar{x}_s$, and $y_s$, hence
of absolute value at most 
$\mbox{max} (|\bar{x}_s|,|y_s|)  \leq \mbox{max}(|x_s|,|x_{s+1}|)= |x_{s+1}|$ 
such that
\[ |f^{(s+1)} (\bar{x}_{s+1})| = 
   \frac{|f^{(s)}(\bar{x}_s) - f^{(s)} (y_s)|}{|\bar{x}_s- y_s|} =
    \frac{|f^{(s)}(\bar{x}_s)|}{|\bar{x}_s- y_s|} \geq
   \frac{1}{2^s} \frac{\prod_{i=s+1}^n |x_i|}{|\bar{x}_s- y_s|}. \]
Since also $|y_s| \leq |x_{s+1}|$, the inequality
$|\bar{x}_s - y_s| \leq 2 |x_{s+1}|$ holds, and consequently,
\[ |f^{(s+1)} (\bar{x}_{s+1})| \geq
   \frac{1}{2^s} \frac{\prod_{i=s+1}^n |x_i|}{2 |x_{s+1}|} =
   \frac{1}{2^{s+1}} \prod_{i=s+2}^n |x_i|, \]
as was to be proved.

\vspace{2mm}

{\sc Step 3:}
If $3n^3 |\bar{x}_s| \leq |x_{2s-1}|$, with $\bar{x}_s$ from step 2, then
\[ \prod_{i=s+1}^n |x_i| \leq 
   2^{s+1} (2s+1) \; \mbox{max}_{s \leq j \leq 3s} |f^{(j)}(0)|. \]

\proof
By Taylor's formula,
\[ f^{(s)} (\bar{x}_s) = 
   \sum_{i=s}^{3s} \frac{f^{(i)}(0)}{(i-s)!} \bar{x}_s^{i-s} +
   \sum_{i=3s+1}^{n} \frac{f^{(i)}(0)}{(i-s)!} \bar{x}_s^{i-s},  \]
hence
\begin{equation} \label{hilf2}
|f^{(s)} (\bar{x}_s)| \leq 
(2s+1) \; \mbox{max}_{s \leq j \leq 3s} |f^{(j)}(0)| +
\sum_{i=3s+1}^{n} \frac{|f^{(i)}(0)|}{(i-s)!} |\bar{x}_s|^{i-s}.  
\end{equation}
Next, for $i \geq 3s+1$, (\ref{hilf0}) implies
\[ |f^{(i)}(0)| |\bar{x}_s|^{i-s} \leq 
   \frac{n!}{(n-i)!} |\bar{x}_s|^{i-s} \prod_{j=i+1}^n |x_j|  =
   \frac{n!}{(n-i)!} \prod_{j=s+1}^n |x_j| 
   \prod_{j=s+1}^{2s} \frac{|\bar{x}_s|}{|x_j|}
   \prod_{j=2s+1}^i \frac{|\bar{x}_s|}{|x_j|}. \]
Since $\bar{x}_s \leq |x_s| \leq |x_j|$ for $j \geq s+1$, and
$\bar{x}_s \leq \frac1{3n^3} |x_{2s+1}| \leq \frac1{3n^3} |x_j|$ for
$j \geq 2s+1$, the above is less or equal
\[ \frac{n!}{(n-i)!} \left(\frac 1{3n^3} \right)^{i-2s} \prod_{j=s+1}^n |x_j| 
   \leq \frac{n^i}{(3n^3)^{i-2s}} \prod_{j=s+1}^n |x_j| \leq
   \left(\frac1 {3}\right)^{i-2s} \prod_{j=s+1}^n |x_j|. \]
Consequently,
\[ \sum_{i=3s+1}^{n} \frac{|f^{(i)}(0)|}{(i-s)!} |\bar{x}_s|^{i-s} \leq
   \prod_{j=s+1}^n |x_j| \sum_{i=3s+1}^n \frac1{(i-s)!}
   \left(\frac1{3}\right)^{i-2s} \leq
   \frac{1}{2 \cdot 3^s} \prod_{j=s+1}^n |x_j|. \]
Together with (\ref{hilf1}), and(\ref{hilf2}), this implies
\[ \frac1{2^s} \prod_{s=1}^n |x_i| \leq 
   (2s+1) \; \mbox{max}_{s \leq j \leq 3s} |f^{(j)}(0)| +
   \frac{1}{2 \cdot 3^s} \prod_{j=s+1}^n |x_j|, \]
hence
\[ \frac1{2^{s+1}} \prod_{i=s+1}^n |x_i| \leq 
   (2s+1) \; \mbox{max}_{s \leq j \leq 3s} |f^{(j)}(0)|, \]
and the claim follows.

\vspace{2mm}

{\sc Step 4:}
For every $\bar{s} \leq s$,
\[ \prod_{i=2s-\bar{s}+1}^{2s} |x_i| \prod_{i=s+1}^n |x_i| \leq \]
\begin{equation} \label{hilf3}
\mbox{max} \left( 4(s+1)(3n^3)^{s+1} \;
         \mbox{max}_{0 \leq j\leq 3s}f^{(j)}(0),
              (3n^3)^{\bar{s}} \prod_{i=s-\bar{s}+1}^n |x_i| \right). 
\end{equation}

\proof
The claim obviously holds for $\bar{s}=0$. Assume (\ref{hilf3}) holds for 
$\bar{s} \leq s-1$. If in (\ref{hilf3}), 
\[ \prod_{i=2s-\bar{s}+1}^{2s} |x_i| \prod_{i=s+1}^n |x_i| \leq
   4(s+1)(3n^3)^{s+1} \;
         \mbox{max}_{0 \leq j\leq 3s}f^{(j)}(0), \]
then
\[ \prod_{i=2s-\bar{s}}^{2s} |x_i| \prod_{i=s+1}^n |x_i| =
   |x_{2s-\bar{s}}| \prod_{i=2s-\bar{s}+1}^{2s} |x_i| \prod_{i=s+1}^n |x_i|
   \leq \]
\[ |x_{2s-\bar{s}}| 4(s+1)(3n^3)^{s+1} \;
         \mbox{max}_{0 \leq j\leq 3s}f^{(j)}(0) \leq
   4(s+1)(3n^3)^{s+1} \;
         \mbox{max}_{0 \leq j\leq 3s}f^{(j)}(0). \]
If in (\ref{hilf3}),
\[ \prod_{i=2s-\bar{s}+1}^{2s} |x_i| \prod_{i=s+1}^n |x_i| \leq
   (3n^3)^{\bar{s}} \prod_{i=s-\bar{s}+1}^n |x_i|, \]
then if  
$3n^3|\bar{x}_{s-\bar{s}}| \leq |x_{2(s-\bar{s})-1}|$, 
step 3, with $s-\bar{s}$ instead of $s$ implies
\[ \prod_{i=s-\bar{s}+1}^n |x_i| \leq 2^{s-\bar{s}+1} (2 (s-\bar{s})+1) 
   \mbox{max}_{s-\bar{s}\leq j\leq 3(s-\bar{s})} |f^{(j)}(0)|, \]
hence,
\[ \prod_{i=2s-\bar{s}}^{2s} |x_i| \prod_{i=s+1}^n |x_i| \leq
   |x_{2s-\bar{s}}| (3n^3)^{\bar{s}} \prod_{i=s-\bar{s}+1}^n |x_i| \leq \]
\[ (3n^3)^{\bar{s}} 2^{s-\bar{s}+1} (2 (s-\bar{s})+1) 
   \mbox{max}_{s-\bar{s}\leq j\leq 3(s-\bar{s})} |f^{(j)}(0)| \leq \]
\[ (3n^3)^{s+1} (2 (s-\bar{s})+1) 
   \mbox{max}_{s-\bar{s}-1\leq j\leq 3(s-\bar{s}-1)} |f^{(j)}(0)|  \leq \]
\[ (3n^3)^{s+1} (2s+1) 
   \mbox{max}_{0 \leq j\leq 3s} |f^{(j)}(0)|, \]
hence (\ref{hilf3}) for $\bar{s}+1$.

If on the other hand 
$3n^3 |\bar{x}_{s-\bar{s}}| \geq |x_{2s-\bar{s}}|$, by (\ref{hilf3}),
either
\[ \prod_{i=2s-\bar{s}}^{2s} |x_i| \prod_{i=s+1}^n |x_i|
   \leq (3n^3)^{\bar{s}+1} \prod_{i=s-\bar{s}}^n |x_i|, \]
or
\[ \prod_{i=2s-\bar{s}}^{2s} |x_i| \prod_{i=s+1}^n |x_i| \leq
   \prod_{i=2s-\bar{s}+1}^{2s} |x_i| \prod_{i=s+1}^n |x_i| \leq \]
\[ 2(s+1)2(3n^3)^{s+1} \;
         \mbox{max}_{0 \leq j \leq 3s}f^{(j)}(0), \]
hence (\ref{hilf3}) for $\bar{s}+1$.

\vspace{2mm}

{\sc Step 5:}
Now, (\ref{hilf3}) for $\bar{s} = s$ reads
\[ \prod_{i=s+1}^{2s} |x_i| \prod_{i=s+1}^n |x_i| \leq
   \mbox{max} \left(2(s+1)(3n^3)^{s+1} \;
         \mbox{max}_{0 \leq j \leq 3s}f^{(j)}(0), (3n^3)^s 
   \prod_{i=1}^n |x_i|\right) = \]
\[ 2(s+1)(3n^3)^{s+1} \;
         \mbox{max}_{0 \leq j \leq 3s}f^{(j)}(0), \]
since $\prod_{i=1}^n |x_i| = f^{(0)}(0)$. Hence,
\[ \prod_{i=s+1}^n |x_i|^2 \leq 
   \prod_{s+1}^{2s} |x_i| \prod_{i=s+1}^n |x_i| \leq
   2(s+1)(3n^3)^{s+1} \;
         \mbox{max}_{0 \leq j \leq 3s}f^{(j)}(0), \]
that is the first inequality of the Lemma.

\proof {\sc of Lemma \ref{hilfzwei}}
Define $x_i := \varphi^{-1}(z_i), i=1, \ldots, n$. Using
\[ f(z) = \prod_{i=1}^n |\theta,z_i| = \frac{|x_i-x|}{\sqrt{1+ |x_i-x|^2}} =
   F(x), \]
the inequality
\begin{equation} \label{hilfb1}
\left|(\partial^I F)(0)\right| \leq n^s \prod_{s+1}^n |\theta, z_i|
\end{equation}
is a straightforward calculation.

Again, we may assume that
$|\theta,z_1| < \cdots < |\theta,z_n|$, and also that $|\theta,z_i|<1$,
for $i=1, \ldots, n$ i.\@ e.\@ the $z_i$ are not at infinity with respect to 
$\theta$. Hence, $|x_1| < \cdots < |x_n|$, and $|x_n| < \infty$.

Since there are only $n$ points $x_i$, there exists a real line
through the origin $L \subset \C^t$ and a permutation $\pi \in \Sigma_n$,
such that with $pr_L$ the projection of $\C^t$ to $L$,
and $y_i = pr_L(x_{\pi i})$,
\[ |y_i| < \cdots < |y_n|, \quad 
   |x_i| \geq |y_{\pi i}| \geq \frac{|x_i|}{n}, \]
and consequently,
\begin{equation} \label{vergl}
|y_i| \leq n |y_{\pi i}| \leq n^2 |y_i|, \quad
|x_i| \leq n |x_{\pi i}| \leq n^2 |x_i|. 
\end{equation}
Let $\partial = \partial_L$ be the directional derivative in the dirction
of $L$, and $\tilde{z_i}$ the point in $\Pe^t(\C)$ corresponding to $y_i$.
Then, with $g(z) = \prod_{i=1}^n |\theta,\tilde{z}_i|$, 
$G = \varphi^* g$, and $s \leq n$,
\begin{equation} \label{um2}
(|\partial_L^s G)(0)| \leq n^s |(\partial_L^s F)(0)|.
\end{equation}

To prove for $s \leq n/3$ that
\[ \prod_{s+1}^n |z_i, \theta| \leq 2 (s+1) (3n^3)^{s+1} n^{n-s}
   \sup_{i\leq 3s} |\partial^i f^*(0)|, \]
assume first that 
$|\theta, z_s| \geq 1/\sqrt{2}$. Since
$\prod_{i=1}^n |\theta, z_i| = F(0)$, one may assume
\[ \prod_{i=s}^n |\theta, z_i| \leq 2 s (3n^3)^s n^{n-s+1}
   \sup_{|I|\leq3(s-1)} \left| (\partial^I F)(0) \right|. \]
and then
\[ \prod_{i=s+1}^n |\theta, z_i| \leq
   \sqrt{2} \prod_{i=s} |\theta, z_i| \leq 
   \sqrt{2} \; 2 s (3n^3)^s n^{n-s+1}
   \sup_{|I|\leq3(s-1)} \left| (\partial^I F)(0) \right| \leq \]
\[ 2 (s+1) (3n^3)^{s+1} n^{n-s}
   \sup_{|I|\leq 3s} \left| (\partial^I f)(0) \right|. \]

Thus, we may from now on assume that $|\theta,z_s| < 1/\sqrt{2}$, 
hence $|y_s| \leq |x_s| < 1$. Then, with $\bar{x}_i := x_1$ if
$|\theta,z_i| < 1/\sqrt{2}$, and 
$\bar{x}_i := x_i/|x_i| \geq |\theta,z_i|$ Lemma \ref{hilf} implies
\[ \prod_{i=s+1}^n |\bar{x}_i| \leq (2s+1) (3n^3)^{s+1}
   \sup_{0 \leq j \leq 3s}  |(\partial^j_L G)(0)|, \]
where $G(x) = \prod_{i=1}^n |\bar{x}_i|$.
Since $|\bar{x}_i| = c_i |\theta,\tilde{z}_i|$, with $c_i \leq 1$ 
and $|\theta,\tilde{z}_i| \leq n|\theta,z_i|$ for 
$i=1, \ldots, n$, we have $|(\partial^j_L G)(0)| \leq |(\partial^j_L F(0)|$.
Further, $\prod_{i=s+1}^n |\theta,z_i| \leq n^n\prod_{i=1}^n |\bar{x}_i|$, 
hence
\[ \prod_{i=s+1}^n |\theta,z_i| \leq (2s+1) (3n^3)^{s+1} n^n
   \sup_{0 \leq j \leq 3s} |(\partial^j_L F)(0)|. \]
Since clearly 
$\sup_{0 \leq j \leq 3s} |(\partial^j_L F)(0)| \leq
\sup_{0 \leq |J| \leq 3s} |(\partial^J F)(0)|$, the Lemma follows.

\proof {\sc of Lemma \ref{hilfdrei}:}
By the definition of the Fubini-Study metric for any multiindex, and
any $i=1, \ldots,n$
\[ |\partial^J \varphi^*(|\theta,\Pe(W)_i|)| \leq 1. \]
The Lemma hence follows from the Leibniz rule.

\end{appendix}


\begin{thebibliography}{0mm}

\bibitem[BGS]{BGS} Bost, Gillet, Soul\'e: Heights of projective varieties and
positive Green forms. JAMS 7,4 (1994)


\bibitem[GS1]{GS1} H.\@ Gillet, C.\@ Soul\'e: Arithmetic intersection theory.
Publications Math.\@ IHES 72 (1987) 243-278

\bibitem[GS2]{GS2} H.\@ Gillet, C.\@ Soul\'e: Characteristic classes for 
algebraic vector bundles with hermitian metric I,II. 
Annals of Mathematics 131 (1990), 163-238.



\bibitem[LR]{LR} M.\@ Laurent, D.\@ Roy: Sur l'approximation alg\'ebrique
en degr\'e de transcendence un.

\bibitem[Ma1]{App1} H.\@ Massold: Diophantine Approximation 
on varieties I: Algebraic Distance and Metric B\'ezout Theorem. 
math.NT/0611715v2


\bibitem[SABK]{SABK} C.\@ Soul\'e, D.\@ Abramovich, J.\@-F.\@ Burnol,
               J.\@ Kramer: Lectures on Arakelov Geometry. 
               Cambridge University Press 1992 

\end{thebibliography}
\end{document}